%% file: main.tex
\documentclass{amsart}

\input{header}

\newcommand{\tempnewpage}{}

\begin{document}


\expandafter\title
{Symmetric and spectral realizations of highly symmetric graphs}
		
\author[M. Winter]{Martin Winter}
\address{Faculty of Mathematics, University of Technology, 09107 Chemnitz, Germany}
\email{martin.winter@mathematik.tu-chemnitz.de\newline\rule{0pt}{1.5cm}\includegraphics[scale=0.7]{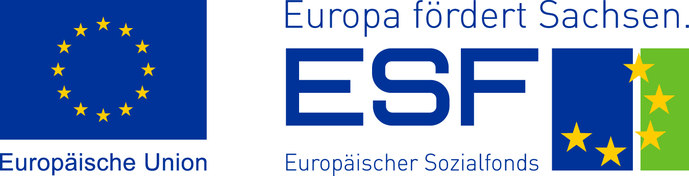}}
	
		
\date{\today}
\begin{abstract}
A \emph{realization} of a graph $G=(V,E)$ is a map $v\: V\to\RR^d$ that assigns to each vertex a point in $d$-dimensional Euclidean space.
We study graph realizations from the perspective of representation theory (expressing certain symmetries), spectral graph theory (satisfying certain self-stress conditions) and rigidity theory (admitting deformations that do not alter the symmetry properties).


We explore the connections between these perspectives, with a focus on realizations of highly symmetric graphs (arc-transitive/distance-transitive) and the question of how much symmetry is necessary to ensure that a realization is balanced, spectral, rigid etc.

%
%
%
%

We include many examples to give a broad overview of the possibilities and restrictions of symmetric and spectral 
graph realizations.
%
\end{abstract}

\maketitle

\input{sec/introduction}

\input{sec/realizations}

\input{sec/future}

{
\par\bigskip
\parindent 0pt
\textbf{Acknowledgements.} The author gratefully acknowledges the support by the funding of the European Union and the Free State of Saxony (ESF).
}


\bibliographystyle{abbrv}
\bibliography{literature}

\newpage
\input{sec/appendix}

\end{document}

%% file: header.tex
\usepackage[T1]{fontenc}
\usepackage[utf8]{inputenc}
\usepackage[USenglish]{babel}
\usepackage{textcase}

\usepackage[
	bookmarks=true,
	plainpages=false,
	linktocpage,
	colorlinks=true,
	citecolor=green!80!black,
	linkcolor=red!70!black,
	filecolor=magenta,
	urlcolor=magenta,
	breaklinks,
	pdfauthor={Martin Winter},
]{hyperref}

\usepackage{amsmath,amsthm,amssymb}
\usepackage{calc, mathtools} 

\usepackage{colortbl,color} 
\usepackage[dvipsnames]{xcolor}
\usepackage{centernot} 
\usepackage{array} 
\usepackage{enumitem,moreenum} 
\usepackage{cite} 
\usepackage[nameinlink,capitalize,noabbrev]{cleveref}
\usepackage{nicefrac}

\usepackage{tikz-cd} 

\usepackage[font=small,labelfont=bf]{caption}

\usepackage{listings}
\newcommand{\ZZ}{\mathbb{Z}}    
\newcommand{\RR}{\mathbb{R}}    
                   
\newcommand{\nozero}{\setminus\{0\}}

\newcommand{\Tsymb}{\top}
\newcommand{\T}{^{\Tsymb}}

\def\^#1{^{(#1)}}
\def\s^#1{^{\smash{(#1)}}}

\newcommand*{\horzbar}{\rule[.5ex]{2.5ex}{0.4pt}}

\newcommand{\wideword}[1]{\quad\text{#1}\quad}
\newcommand{\wideand}{\wideword{and}}

\def\:{\colon}

\newcommand{\cupdot}{\mathbin{\mathaccent\cdot\cup}}

\newcommand{\mylabel}{\upshape(\textit{\roman*})}
\newenvironment{myenumerate}{\begin{enumerate}[label=\mylabel]}{\end{enumerate}}

\newcommand{\freespace}{\kern.07em} 
\newcommand{\free}{\freespace\cdot\freespace} 

\newcommand{\card}[1]{|#1|}

\newcommand{\enquote}[1]{``#1''}                                 

\newcommand{\ul}[1]{\underline{\smash{#1}}}

\newcommand{\msays}[1]{{\footnotesize\textcolor{red}{#1}}}
\newcommand{\TODO}{\msays{TODO}}

\theoremstyle{plain}  
\newtheorem{theorem}{Theorem}[section]
\newtheorem{corollary}[theorem]{Corollary}
\newtheorem{lemma}[theorem]{Lemma}
\newtheorem{proposition}[theorem]{Proposition}

\theoremstyle{definition} 
\newtheorem{definition}[theorem]{Definition}
\newtheorem{example}[theorem]{Example}
\newtheorem{remark}[theorem]{Remark}
\newtheorem{question}[theorem]{Question}
\newtheorem{observation}[theorem]{Observation}

\crefname{theorem}{Theorem}{Theorems}
\crefname{proposition}{Proposition}{Propositions}
\crefname{lemma}{Lemma}{Lemmas}
\crefname{corollary}{Corollary}{Corollaries}
\crefname{remark}{Remark}{Remarks}
\crefname{example}{Example}{Examples}
\crefname{definition}{Definition}{Definitions}
\crefname{problem}{Problem}{Problems}
\crefname{observation}{Observation}{Observation}
\crefname{construction}{Construction}{Construction}

\DeclareMathOperator{\aff}{aff}

\DeclareMathOperator{\Aut}{Aut}
\DeclareMathOperator{\diam}{diam}

\DeclareMathOperator{\Ortho}{O}

\DeclareMathOperator{\GL}{GL}

\DeclareMathOperator{\Span}{span}

\DeclareMathOperator{\rank}{rank}
\DeclareMathOperator{\tr}{tr}

\DeclareMathOperator{\Id}{Id}
\DeclareMathOperator{\id}{id}

\DeclareMathOperator{\Fix}{Fix}
\DeclareMathOperator{\dist}{dist}
\DeclareMathOperator{\Eig}{Eig}
\DeclareMathOperator{\Spec}{Spec}
\DeclareMathOperator{\Sym}{Sym}
\DeclareMathOperator{\Perm}{Perm}   	

\let\<=\langle
\let\>=\rangle
\let\x=\times

\def\...{...}
\newcommand{\shortStyle}{\textit}
\newcommand{\ie}{\shortStyle{i.e.,}}

\newcommand{\eg}{\shortStyle{e.g.}}

\newcommand{\wrt}{\shortStyle{w.r.t.}}

\newcommand{\cf}{\shortStyle{cf.}}

\newcommand{\resp}{resp.}

\let\angle=\measuredangle

\makeatletter
\renewcommand*{\eqref}[1]{%
  \hyperref[{#1}]{\textup{\tagform@{\ref*{#1}}}}%
}
\makeatother

\numberwithin{equation}{section}

%% file: sec/introduction.tex
\section{Introduction}

Throughout the paper, $G=(V,E)$ denotes a (simple, undirected) graph with~vertex set $V=\{1,...,n\}$. 
In general, we assume that $G$ has many symmetries, \ie~has a large symmetry group $\Aut(G)\subseteq\Sym(V)$.

A $d$-dimensional \emph{(graph) realization} \mbox{$v\:V\to\RR^d$} is a (not necessarily injective) embedding of the vertices of $G$ into the $d$-dimensional Euclidean space, and one~can consider this as an embedding of the whole graph by imagining the edges embedded as straight lines between the vertices: 
\begin{figure}[!h]
\begin{center}
\includegraphics[width=0.8\textwidth]{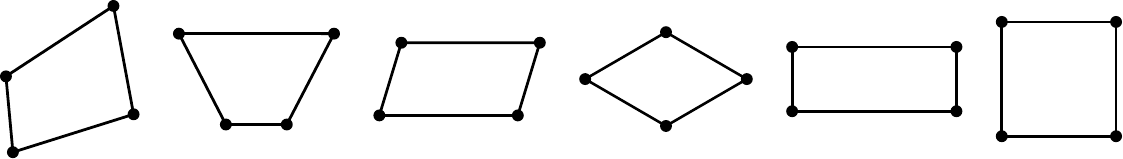}
\end{center}
\caption{Several 2-dimensional realizations of the 4-cycle graphs $G=C_4$.}
\label{fig:C4}
\end{figure}

In this paper we shall discuss 
various classes of such \mbox{realizations}, based on ideas from representation theory (\emph{symmetric} and \emph{rigid} realizations), and~spec\-tral graph theory (\emph{balanced} and \emph{spectral} realizations).
For a rule of thumb: people are interested in symmetric realizations (because they tell a lot about the structure of $G$), and we can use spectral realization to obtain such (because they are fast and easily computed).
 
That spectral realizations of highly symmetric graphs are always highly symmetric (in a way we make precise below) is well-understood and frequently utilized (\eg\ in graph-drawing algorithms).
It is the other direction which provides interesting unanswered questions: \enquote{How much symmetry is necessary for a realization to be spectral?}.
We are trying to answer this question for several classes of highly symmetric graphs, most successfully in the case of distance-transitive graphs.

Our overall goal is to explore and to better understand the general \mbox{connections} between spectral and symmetric realizations. 
These connections were also \mbox{studied} in \cite{berkolaiko2018eigenspaces} for quantum and edge-weighted graphs, as well as in \cite{du2013graph}. In Section 4 of \cite{du2013graph} the author remarked that the \emph{Petersen graph} enjoys a remarkable property, namely, that all its eigenspaces are $\Aut(G)$-irreducible. This is part of a larger pattern that we discuss in \cref{sec:distance_transitive}.
Another classic reference for spectral properties of highly symmetric graphs is \cite{lovasz1975spectra}.
Our investigation is supplemented with numerous examples that hopefully provide the reader with some intuition and visualization for these relations. 


A secondary goal is to demonstrate that the concept of the \emph{arrangement space} (as intro\-duced in \cite{winter2019geometry}, see also \cref{def:arrangement_space}) is the appropriate tool for defining and relating these diffe\-rent concepts within a \mbox{common} terminology.~Most~of our results have a second inter\-pretation in terms of arrangement spaces, and many proofs make use of this language.

To better explain our findings, we briefly introduce the most relevant terminology.




\subsection{Symmetric realizations}

Informally, a \enquote{symmetric realization} is a realization that
manifests all (or many) of the combinatorial symmetries of $G$ as geometric symmetries. 
Formally, we have 

\begin{definition}
\label{def:symmetric_realization}
For some group $\Sigma\subseteq\Aut(G)$ of symmetries of $G$, a realization is called \emph{$\Sigma$-realization} if there is a linear (orthogonal) representation\footnote{Check out \cref{sec:appendix_rep_theory} for the relevant basics of representation theory.} $T\:\Sigma\to\Ortho(\RR^d)$ with
\begin{equation}\label{eq:symmetric}
T_\sigma v_i=v_{\sigma(i)},\quad\text{for all $i\in V$ and $\sigma\in\Sigma$}.
\end{equation}
%
\end{definition}

Similar constructions, though with a distinct philosophy, are known in finite frame theory as \emph{group frames} (see \cite{waldron2018introduction} for an introduction in finite frame theory, and \cite{waldron2018group} for group frames in specific).


A realization can be at most as symmetric as its underlying graph, which happens if we have $\Sigma=\Aut(G)$.
But by choosing appropriate subgroups $\Sigma\subseteq\Aut(G)$, we have fine control over the kind of symmetries that we want to realize.
At least most of the time.
Occasionally, any $\Sigma$-realization is also a $\Sigma'$-realization for some larger group $\Sigma'\supset\Sigma$.
We say that $\Sigma$ cannot be geometrically realized \enquote{in isolation}.
We address this phenomenon in the later sections.

In the example in \cref{fig:C4}, the left-most realizations is void of any symmetries, and can be a $\Sigma$-realization only for $\Sigma=\{\id\}$.
The right-most realization is as symmetric as possible, or, as we are going to say, is an $\Aut(G)$-realization.
In between, we also find realizations that only realize either vertex- or edge-transitivity.
Apparently, in this case, a separation of certain sub-symmetries is possible.

Among the symmetric realizations we further distinguish the \emph{irreducible} realiza\-tions (those, for which $T$ is an irreducible representation),
and \emph{rigid} realizations (those, which cannot be continuously deformed without immediately becoming less symmetric; see \cref{def:rigid}).
This notion of rigidity has to be distinguished from the one studies in rigidity theory of frameworks. We do not necessarily require that the edge-length are fixed during a deformation (but this is implicitly true if \eg\ $\Sigma$ acts edge-transitively on $G$; see \cref{sec:arc_transitive}).
Still, frameworks with symmetry constraints have been investigated before \cite{schulze2010symmetry,malestein2014generic,nixon2016symmetry}.



\subsection{Balanced and spectral realizations}

Balanced and spectral realizations on the other hand are realizations related to the spectral properties of $G$ (\ie~the eigenvalues and eigenvectors of its adjacency matrix $A$)\footnote{Check out \cref{sec:appendix_spec_theory} for the relevant basics of spectral graph theory.} and do not impose any symmetry constraints a priori.

Spectral realizations have been around for quite some time and were utilized in diverse contexts.
Tracing the historical roots of this notion is beyond the scope of this article.
Naming only a few, there have been applications in data visualization (in particular, graph drawings \cite{koren2003spectral}), semi-definite optimization (\eg\ eigenvalue optimization \cite{boyd2004fastest,goring2008embedded}), geometric combinatorics (\eg\ for equiangular lines \cite{lemmens1991equiangular} and balanced point arrangements \cite{cohn2010point}) as well as polytope theory (in the form of \emph{eigenpolytopes} \cite{godsil1978graphs}).
They are further related to the Lovász theta function \cite{lovasz1979shannon} and the Colin de Verdière graph invariant \cite{van1999colin}.

There are serveral ways to introduce spectral realizations, one of which is via the balanced realizations:

\begin{definition}
\label{def:balanced}
A realization $v\: V\to\RR^d$ is said to be \emph{$\theta$-balanced} (or just \emph{balanced}) for some $\theta\in\RR$ if
\begin{equation}
\label{eq:balanced}
\sum_{\mathclap{j\in N(i)}} v_j= \theta v_i,\quad\text{for all $i\in V$}.
\end{equation}
\end{definition}

Equation \eqref{eq:balanced} can be interpreted as a self-stress condition (from which the name \enquote{balanced} can be motived). 
Alternatively, and this is the perspective of this paper, \eqref{eq:balanced} can be interpreted as an eigenvalue equation for the adjacency matrix of $G$.
More precisely, if we define the so-called \emph{arrangement matrix}
\begin{equation}\label{eq:arrangement_matrix}
M :=\begin{pmatrix}
	\;\horzbar\!\!\!\! & v_1\T & \!\!\!\!\horzbar\;\; \\
	& \vdots & \\[0.4ex]
	\;\horzbar\!\!\!\! & v_n\T & \!\!\!\!\horzbar\;\;
\end{pmatrix}\mathrlap{\in\RR^{n\times d},}
\end{equation}
in which the $v_i$ are the rows, then \eqref{eq:balanced} reads $AM=\theta M$.
So, $\theta$ is an eigenvalue of the adjacency matrix $A$, and the columns of $M$ are (some) \mbox{corresponding~eigen}\-vectors.
In the extreme case, when the columns of $M$ are a complete set of  $\theta$-eigen\-vectors, \ie\ they span the $\theta$-eigenspace $\Eig_G(\theta)$, we speak of a \emph{spectral realization}.

%
%
\begin{definition}
\label{def:spectral}
A $\theta$-balanced realization $v\:V\to\RR^d$ is called \emph{$\theta$-spectral} (or just \emph{spectral}) if the multiplicity of $\theta\in\Spec(G)$ is $d$.
\end{definition}

The previous discussion already gives a description on how to construct $\theta$-spectral realizations: find a basis of the $\theta$-eigenspace, put them as columns of $M$, and read of the $v_i$ in the rows.
The reader can also find attached a short Mathematica script for computing spectral realization in \cref{sec:appendix_mathematica}.

We also mentioned before that spectral realizations are as symmetric as possible, which we can now state as \enquote{spectral realizations are $\Aut(G)$-realizations} (which is well known, but we include a proof in \cref{res:spectral_then_symmetric}).

\subsection{Outline of the paper}

In \cref{sec:realizations} we are setting the stage for our investigations.
We define the notion of the \emph{arrangement space} of a realization and discuss how it can be used to study graph realizations up to orthogonal transformations (see \cref{res:equivalence}), and how it characterizes symmetric, balanced and spectral realizations (see \cref{res:arrangement_space_symmetric} and \cref{res:arrangement_space_spectral}).
We prove briefly that spectral realizations are always as symmetric as the underlying graph (see \cref{res:spectral_then_symmetric}).
We introduce the notions of \emph{deformations} and \emph{rigidity} for realizations, and explain how these relate to spectral and representation theoretic properties of $G$ and $\Aut(G)$.
For example, we prove that every irreducible realization can be continuously deformed into a balanced realization (see \cref{res:deform_into_balanced}).
This section contains numerous references to a previous paper \cite{winter2019geometry} in which related notions where discussed for point arrangements (instead of graph realizations). Several proofs are cited from this source.

The rest of the papers now investigates the changes to the previous results if we impose stronger and stronger symmetry restrictions on out realizations.

In \cref{sec:vertex_transitive} we explore the consequences of vertex-transitivity, or the lack thereof. 
We argue that questions about rigidity are almost always only meaningful if we assume vertex-transitivity (\cref{res:vertex_transitive_rigid} and \cref{res:vertex_transitive_rigid_cor}).

In \cref{sec:arc_transitive} we investigate edge- and arc-transitive realizations.
\mbox{In the latter~case,} many metric properties of the realization can already be determined from $G$ and its eigenvalues (see \cref{res:cosine_and_length}).
General (irreducible) arc-transitive realizations are still not too well-behaved,~\eg\ they are not necessarily rigid or spectral (and we do not know whether they are necessarily balanced).
We can show that under mild assumption (namely, \emph{full local dimension}, see \cref{def:full_local_dimension}) an arc-transitive realization is rigid, irreducible and balanced (see \cref{res:full_local_dimension_is_good}).

The final section, \cref{sec:distance_transitive}, is devoted to realizations with exceptionally strong symmetry, namely, \emph{distance-transitivity}.
We explain how the distance-transitive realizations of a graph are completely accessible via spectral realizations (see \cref{res:irreducibe_distance_transitive_is_spectral}).
\cref{res:distance_transitive_consequences} shows that all distance-transitive realizations are rigid and realize all the symmetries of the underlying graph (or in other words, distance-transitivity cannot be realized in isolation).
We explore a generalization of this result beyond distance-transitivity with the help of \emph{cosine vectors} (see \cref{res:cosine_vectors}).

%% file: sec/realizations.tex
\section{Realizations and arrangement spaces}
\label{sec:realizations}

From this section on, let $v\: V\to\RR^d$ denote a realization of \emph{full dimension}, which means that $\rank v$ $:=\dim \Span\{v_1,...,v_n\}=d$.
Equivalently, the arrangement~matrix $M$ (as defined in \eqref{eq:arrangement_matrix}) has rank $d$.

\subsection{Normalized and spherical realizations}

We do care about metric properties of realizations, as lengths and angles, but we do not care about the exact positioning of each vertex in space, that is, we do not care about the orientation of the realization.
We now introduce some terminology that enables us to study realizations up to orthogonal transformations.


\begin{definition}\label{def:arrangement_space}
Given a realization $v$ with arrangement matrix $M$,
\begin{myenumerate}
	\item if $M\T M=\alpha\Id$ for some $\alpha>0$, then $v$ is called \emph{spherical}.
	\item if $M\T M=\Id$, then $v$ is called \emph{normalized}.
	\item the column span $U:=\Span M\subseteq\RR^n$ is called \emph{arrangement space} of $v$.
	\item realizations with the same arrangement space are called \emph{equivalent}.
\end{myenumerate}
\end{definition}

\begin{figure}[h!]
\centering
\includegraphics[width=0.43\textwidth]{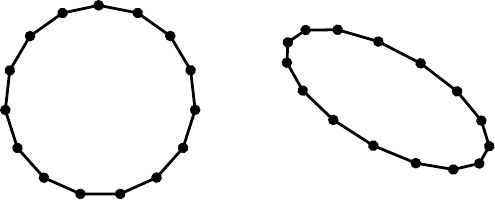}
\caption{Two realizations of a cycle graph, one of which is \enquote{spherical} (left), and one of which is not (right).}
\label{fig:spherical}
\end{figure}

Spherical (and normalized) realizations are always of full dimension, since $\rank M$ $=\rank (M\T\! M)=\rank(\alpha\Id)=d$.
The reason for the use of the word \enquote{spherical} is illustrated in \cref{fig:spherical}.

\noindent
Normalized realizations additionally prescribe a certain scale.

Two (full-dimensional) realizations of the same graph are equivalent (\ie\ have the same arrangement space) if and only if they are related by an invertible linear transformation. This follows from well-known facts in linear algebra: two matrices $M,\bar M\in\RR^{n\x d}$ have the same column span if and only if $M=\bar MT$ for some $T\in\GL(\RR^d)$.
Thus, working with only the arrangement space of a realization already provides a tool for considering realizations up to invertible linear transformations.

But we actually want \enquote{up to orthogonal transformations}.
We can achieve this~via normalized realizations:

%
%
\begin{theorem}[\!\cite{winter2019geometry}, Theorem 3.2]
\label{res:equivalence}
Two normalized realizations $v,w\: V\to\RR^d$ of $G$ are equivalent (\ie\ have the same arrangement space) if and only if they are related by an orthogonal transformation $T\in\Ortho(\RR^d)$, that is, $v_i=T w_i$ for all $i\in V$. 
\end{theorem}

Motivated by \cref{res:equivalence}, we mostly restrict to study normalized realizations. This comes with no major loss of generality.
Foremost, every realization is equivalent to a normalized realization.
If a realization has arrangement space $U\subseteq\RR^n$, we can define a second realization $\bar v$ whose arrangement matrix $\bar M\in\RR^{n\x d}$ has as columns an orthonormal basis of $U$.
Then $v$ and $\bar v$ are equivalent (both have $U$ as arrangement space), and since $\bar M\T\! \bar M=\Id$, we find that $\bar v$ is normalized.

One might object that some realization that are interesting for their symmetry are not normalized (\eg\ the rhombus realization in \cref{fig:C4}).
Later on, \cref{res:irreducible_spherical} gives a reason not to worry about this.


%
%


\begin{remark}
As a consequence of \cref{res:equivalence}, metric properties of normalized realizations are uniquely determined by the arrangement space: consider \eg\ the \emph{radius} $r(v)$ defined as follows:
\begin{equation}\label{eq:radius}
[r(v)]^2 := \frac 1n\sum_{\mathclap{i\in V}} \|v_i\|^2 = \frac1n \tr\underbrace{(M\T\! M)}_{=\Id} =\frac{\dim U} n = \frac dn.
\end{equation}
If all vertices of the normalized realization are on a common sphere around the origin (\eg\ as they are for every vertex-transitive realization, see \cref{res:vertex_transitive} in \cref{sec:vertex_transitive}), then the radius of this sphere is given by $r(v)=d/n$.
\end{remark}

\subsection{Symmetric realizations}

In the context of symmetric realizations there is another reason not to worry about a restriction to normalized (or spherical) realizations:
%
%
\begin{theorem}[\!\!\cite{winter2019geometry}, Proposition 4.7]
\label{res:irreducible_spherical}
An irreducible $\Sigma$-realization is spherical.
\end{theorem}
%

\begin{remark}
\label{rem:irreducible_decomposition}
The irreducible $\Sigma$-realizations are the building blocks of general~$\Sigma$-realizations: for a $\Sigma$-realization $v$ with representation $T$, the ambient space decomposes as $\RR^d=W_1\oplus\cdots\oplus W_m$ into a direct sum of pairwise orthogonal $T$-irreducible subspaces $W_k\subseteq\RR^d$. 
The $k$-th irreducible component $v^{(k)}\:V\to\RR^{d_k}$ of $v$ ($d_k$ being the dimension of $W_k$) then is (equivalent to) the orthogonal projection of $v$ onto the subspace $W_k$.

If we assume that $v$ was conveniently oriented, so that the $W_k$ are contained in the coordinate planes, then we can write $v$ as a concatenation 
\begin{equation}\label{eq:irreducible_decomposition}
v_i=\big(v_i^{(1)},...,v_i^{(m)}\big)\in\RR^{d_1+\cdots+ d_m}.
\end{equation}
%
%
If the $\Sigma$-representation of $v^{\smash{(k)}}$ is $T^{\smash{(k)}}\:\Sigma\to\Ortho(\RR^{d_k})$, then the representation $T$ of $v$ can be written in block form
\begin{equation}
\label{eq:representation_decomposition}
T_\sigma = \begin{pmatrix}
T_\sigma^{(1)} & & \\
 & \ddots & \\
 & & T_\sigma^{(m)}
\end{pmatrix}\in\Ortho(\RR^{d}),\quad\text{for all $\sigma\in\Sigma$}.
\end{equation}

Further discussion on reducible and irreducible realization (in the form of point arrangements) can be found in \cite{winter2019geometry} (for example, see Lemma 4.6 in \cite{winter2019geometry} for a discussion of several statements equivalent to being reducible).
\end{remark}

%
%





Now, the second big use of the arrangement space is that it provides an \mbox{alternative} point of view on many properties of realizations.
For example, symmetric realizations are characterized as follows: 

%

\begin{theorem}
\label{res:arrangement_space_symmetric}
Given a (full-dimensional) realization $v\: V\to\RR^d$ with arrangement space $U\subseteq\RR^n$, then
\quad
\begin{myenumerate}
	\item if $v$ is a $\Sigma$-realization, then $U$ is a $\Sigma$-invariant subspace of $\RR^n$ (\!\emph{\cite{winter2019geometry}}, Theorem 4.8). If $v$ is spherical, then the converse holds (\!\emph{\cite{winter2019geometry}}, Theorem 4.9).
	\item if $v$ is a $\Sigma$-realization, then it is irreducible if and only if $U$ is $\Sigma$-irreducible as invariant subspace of $\RR^n$. \!(\!\emph{\cite{winter2019geometry}}, Theorem 4.11)
\end{myenumerate}
\end{theorem}

The condition of being spherical in the converse of \cref{res:arrangement_space_symmetric} $(i)$ is necessary:~a rectangle, a rhombus and a square all have the same arrangement space (they are all linear transformations of each other).
This arrangement space is $\Aut(G)$-invariant for $G=C_4$ (the cycle graph on four vertices), but only the square is an $\Aut(G)$-realization of $G$.
It is also the only realization of these which is spherical.

\subsection{Balanced and spectral realizations}

Likewise, balanced and spectral realizations are characterized using the arrange\-ment space as follows:
\begin{theorem}\label{res:arrangement_space_spectral}
Given a (full-dimensional) realization $v\: V\to\RR^d$ with arrangement space \mbox{$U\subseteq\RR^n$,~then}
\quad
\begin{myenumerate}
	\item $v$ is $\theta$-balanced if and only if $U\subseteq\Eig_G(\theta)$ (where $\Eig_G(\theta)$ denotes the~$\theta$-eigenspace of $G$).
	\item $v$ is $\theta$-spectral (that is, $v$ is a $\theta$-realization) if and only~if~$U=$ $\Eig_G(\theta)$.
\end{myenumerate}
\begin{proof}
The defining equality \eqref{eq:balanced} for being balanced can be written as $AM=\theta M$ ($A$ being the the adjacency matrix of $G$, and $M$ the arrangement matrix of~$v$).~In this form it is clear that $\theta$ is an eigenvalue of $A$, and that the columns of $M$ are corresponding eigenvectors.
Since the arrangement space $U$ is the column span of $M$,~equation~\eqref{eq:balanced} is equivalent to $U\subseteq\Eig_G(\theta)$. This proves $(i)$.

The dimension of $v$ (assuming full dimension) equals the rank of $M$, which equals the dimension of $U$.
Thus, if the dimension of a balanced realization $v$ agrees with the dimension of the eigenspace, we must have $U=\Eig_G(\theta)$. This proves $(ii)$.
\end{proof}
\end{theorem}

\Cref{res:arrangement_space_spectral} justifies that we speak of \emph{the} $\theta$-realization of $G$, as any two such~rea\-lizations have the same arrangement space, hence (if normalized) differ only by an orthogonal transformation (by \cref{res:equivalence}).

In general, if not mentioned other wise, speaking of \emph{the} $\theta$-realization, we mean the uniquely determined (up to orientation) normalized realization with $U=\Eig_G(\theta)$.
Also, balanced and spectral realizations will be assumed to be at least spherical.

\begin{observation}\label{res:eigenvectors_from_realization}
If, somehow, we got our hands on a $\theta$-balanced realization $v$, it is straight forward to extract $\theta$-eigenvectors of $G$ from that: the elements of the arrangement space $U\subseteq\RR^n$ are exactly the~vectors \mbox{$u\in\RR^n$} with components
\begin{equation}\label{eq:eigenvector}
u_i=\<x,v_i\>,\quad\text{for all $i\in V$},
\end{equation}
for some $x\in\RR^d$. This is because equation \eqref{eq:eigenvector} is equivalent to $u=Mx$ (where~$M$ is the arrangement matrix of $v$).
Therefore $u\in\Span M=U$, and since $v$ is balanced we have $U\subseteq \Eig_G(\theta)$, which shows that $u$ is a $\theta$-eigenvector.
\end{observation}

In the language of arrangement spaces it follows immediately that spectral~reali\-zations have all the symmetries of $G$, that is, are $\Aut(G)$-realizations.
This follows from the well-known fact that eigenspaces of $G$ are $\Aut(G)$-inva\-riant: 

\begin{corollary}\label{res:spectral_then_symmetric}
A (spherical) spectral realization $v$ is an $\Aut(G)$-realization.
\begin{proof}
Let $v$ be a (spherical) spectral realization of $G$.

By \cref{res:arrangement_space_spectral} $(ii)$, the arrangement space $U\subseteq\RR^n$ of $v$ is an eigenspace, and by this, $\Aut(G)$-invariant (this is well-known, but we included an argument below).
By \cref{res:arrangement_space_symmetric} $(i)$ (and since $v$ is spherical) $v$ is thus an $\Aut(G)$-realization.

To see that $U$ is $\Aut(G)$-invariant, recall that the combinatorial symmetries of $G$ are characterized by $A\Pi_\sigma=\Pi_\sigma A$ for all $\sigma\in\Sym(V)$ (where $\Pi_\sigma\in\Perm(n)$~denotes the permutation matrix associated with the permutation $\sigma$).
For every $u\in U=\Eig_G(\theta)$ and $\sigma\in\Aut(G)$ then holds
$$A(\Pi_\sigma u)=(A\Pi_\sigma) u = (\Pi_\sigma A) u = \Pi_\sigma (Au)= \Pi_\sigma (\theta u) = \theta(\Pi_\sigma u),$$
thus $\Pi_\sigma u\in\Eig_G(\theta)=U$ and $U$ is $\Aut(G)$-invariant.
\end{proof}
\end{corollary}


This can be applied in practice: there exist fast an robust algorithms for computing eigenvalues and eigenvectors, and so spectral methods are often the tool~of~choice for obtaining symmetric realizations of a graph (for example, in graph drawings algorithms, discussed \eg\ in \cite{koren2005drawing}). 
Explicitly, the steps are as follows: one computes an \emph{orthonormal basis} $\{u_1,...,u_d\}\subseteq\RR^n$ of the $\theta$-eigenspace of $G$ and defines the matrix $M:=(u_1,...,u_d)\in\RR^{d\x n}$ in which the $u_i$ are the columns.
This is the arrangement matrix of the $\theta$-realization of $G$, and as such, an $\Aut(G)$-realization (by \mbox{\cref{res:spectral_then_symmetric}})\footnote{See \cref{sec:appendix_mathematica} for an implementation in Mathematica.
}.
This technique alone cannot control the dimension of the realization (it is always the multiplicity of $\theta$), and so some adjustments might be necessary depending on the graph and the setting.

\begin{remark}\label{rem:theta_1}
The largest eigenvalue of a graph is always of multiplicity one (see \cref{sec:appendix_spec_theory}).
Therefore, the corresponding $\theta_1$-realization is always of dimension $d=1$, which is not very interesting.
One therefore considers the $\theta_2$-realizations as the first interesting realizations of a graph.
\end{remark}

\begin{example}\label{ex:dodecahedron}
The spectrum of the edge-graph of the \emph{dodecahedron} is 
$$\{(-\sqrt5)^3,(-2)^4,0^4,1^5,(\sqrt 5)^3,3^1\},\quad\text{(exponents denote multiplicites)}.$$ 
Note in particular the two eigenvalues of multiplicity three, $\sqrt 5$ and $-\sqrt 5$.
The both corresponding 3-dimensional spectral realizations are shown below:

\begin{center}
\includegraphics[width=0.55\textwidth]{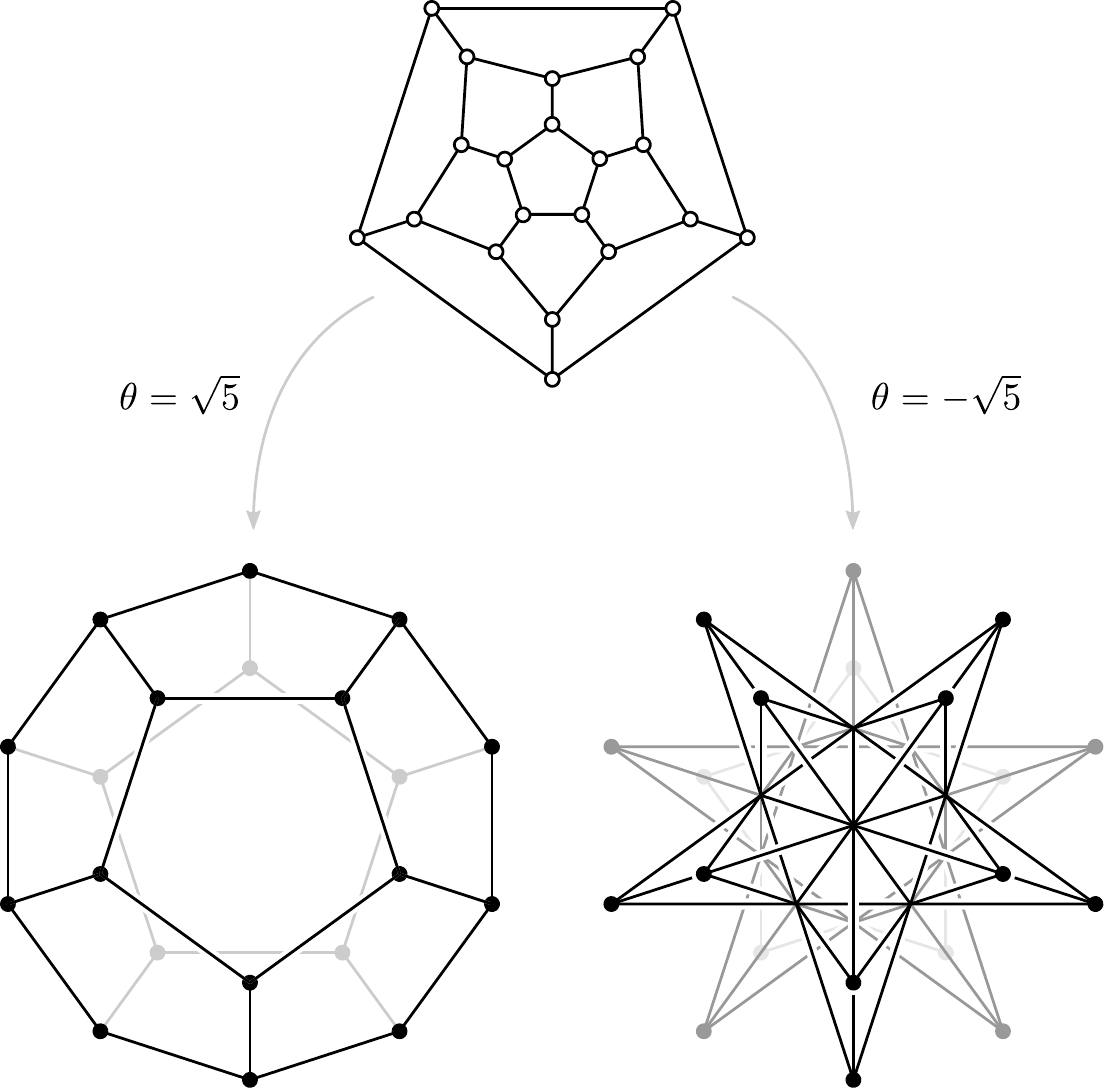}
\end{center}
We observe that the $\sqrt5$-realization gives exactly the skeleton of the regular dodecahedron. 
Note that this is the realization to the \emph{second largest} eigenvalue $\theta_2=\sqrt 5$ of the graph.
This is not completely unexpected.
There are heuristic arguments (\eg\ via nodal domain) that suggest that if the skeleton of a polytope is a spectral realization, then it is the $\theta_2$-realization.
However, a rigorous proof of this observation is still missing.


More evidence for the \enquote{specialness} of $\theta_2$ is provided by the observation that the same phenomenon occurs for the edge-graphs of all the regular polytopes.
%
This was shown by Licata and Powers \cite{licata1986surprising} for all regular polytopes excluding the exceptional 4-dimensional regular polytopes (the 24-cell, 120-cell and 600-cell).
We later close this gap with a note in \cref{ex:24_cell}, or via an alternative approach in \cref{ex:24_cell_cuboctahedron}.
In general however it is an open question which polytope skeleta~can be obtained as spectral realizations.
\end{example}

Besides the fact that we cannot directly control the dimension of the realization, there~are other drawbacks in using spectral realizations when our main focus is on symmetric realizations.
Most notably, an eigenspace, while $\Aut(G)$-invariant, might not be $\Aut(G)$-irreducible. 
If this happens, then the actual irreducible $\Aut(G)$-realizations are inaccessible by spectral method alone (this can probably be fixed, which is briefly discussed in \cref{sec:orbitals}, but this involves additional techniques from computational group theory).
In the language of realizations we can say that some irreducible realizations are not spectral, but only balanced.
And it can be worse: an irreducible symmetric realization must not even be balanced.
We provide examples for both situations.

\begin{example}\label{ex:hexagonal_prism}
Let $G$ be the edge-graph of the hexagonal prism with vertex~set $V=\{1,...,6,1',...,6'\}$,
labels assigned to the vertices as shown in \cref{fig:prism} (left). 
\begin{figure}[!h]
\begin{center}
\includegraphics[width=0.65\textwidth]{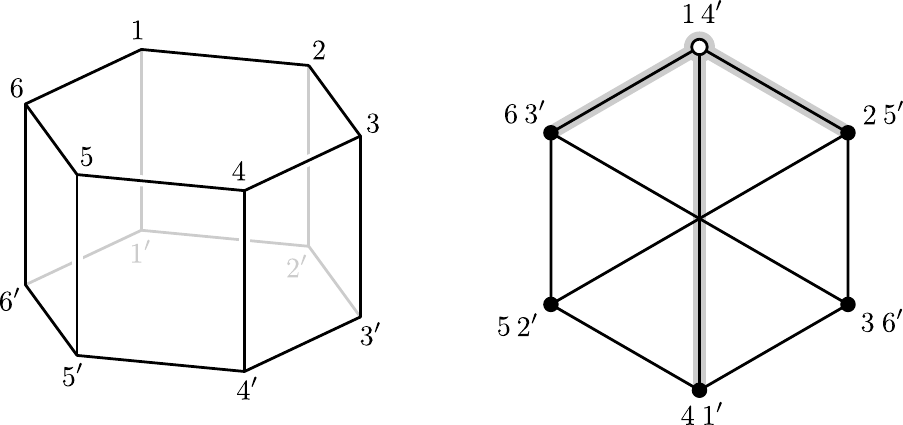}
\end{center}
\caption{The edge-graph of the prism, and an irreducible non-spectral 2-dimensional $\Aut(G)$-realization to the eigenvalue 0.}
\label{fig:prism}
\end{figure}

The spectrum of $G$ is $\{(-3)^1,(-2)^2,(-1)^1,0^4,1^1,2^2,3^1\}$.
Note in \mbox{particular} the  eigenvalue 0 of multiplicity \emph{four}.
We show that the 0-eigenspace is not~irreducible, by constructing a balanced $\Aut(G)$-realization to this eigenvalue of dimension less than four.

Consider the realization in \cref{fig:prism} (right).
The un-dashed vertices are placed in the shape of a hexagon centered at the origin, and each dashed vertex is placed opposite to its un-dashed neighbor, that is, $v_{i'}=-v_i$ (the figure shows only six points because the vertices are mapped on top of each other; the highlighted lines are the images of the edges incident to the vertex $1\in V$).

One checks that this is an $\Aut(G)$-realization of $G$.
One checks further, that this realization is balanced with eigenvalue zero: the three neighbors of each vertex span a regular triangle whose barycenter is the origin.
%
It cannot be a spectral realization since its dimension is not four.
\end{example}

%

%
%

\begin{example}\label{ex:truncated_tetrahedron}
The \emph{truncated tetrahedron} is a polyhedron obtained from the regular tetrahedron by cutting of each vertex.
The image below shows serveral distinct realizations of this polyhedron:
\begin{center}
\includegraphics[width=0.6\textwidth]{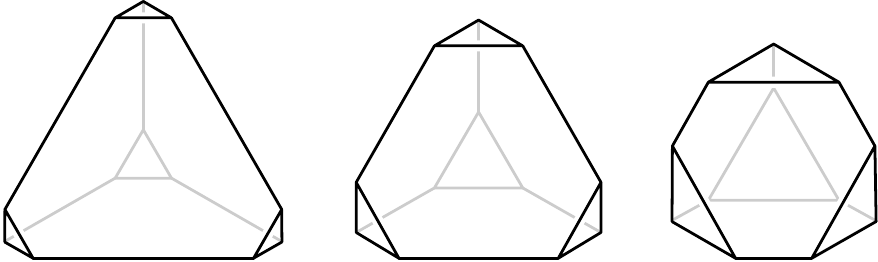}
\end{center}

All of these realizations have the full symmetry of the tetrahedron, which implies that they are $\Aut(G)$-realization of the edge-graph $G$.

Now, this degree of freedom gives an uncountable family of mutually non-equivalent $\Aut(G)$-realizations, each of which corresponds to a distinct $\Aut(G)$-invariant subspace of $\RR^n$ (by \cref{res:arrangement_space_symmetric}).
The spectrum of $G$ is $\{(-2)^3,(-1)^3,0^2,2^3,3^1\}$,~and we see that there are exactly three 3-dimensional balanced realizations of $G$ (which must then also be spectral).
%
%
Consequently, almost all of the previously discussed uncountably many $\Aut(G)$-reali\-zations are not balanced.
%

%
%
\end{example}

This last example contains an instance of a continuous transition between non-equivalent realizations.

\subsection{Deformations and rigidity}

Let $\mathcal R_d(G,\Sigma)$ be the space of all full-dimensional $\Sigma$-realizations of $G$ of dimension $d$. Since realizations can be naturally associated with their arrangement matrices, we can consider $\mathcal R_d(G,\Sigma)$ as a subspace of $\RR^{n\x d}$ equipped with the  subspace topo\-logy.

\begin{definition}\label{def:rigid}
\quad 
\begin{myenumerate}
	\item Given two realizations $v,w\in\mathcal R_d(G,\Sigma)$, we say that these realizations can be \emph{deformed} into each other, if there exists a continuous curve $v(\free)\:[0,1]\to\mathcal R_d(G,\Sigma)$ with $v(0)=v$ and $v(1)=w$.
	The curve $v(t)$ is then called a~\emph{defor\-ma\-tion} between $v$ and $w$.
	\item A realization $v\in\mathcal R_d(G,\Sigma)$ is called \emph{rigid} if it cannot be deformed into~a~non-equivalent realization. It~is~called \emph{flexible} otherwise.
\end{myenumerate}
\end{definition}

Deformations and rigidity of realizations (in the form of symmetric point arrange\-ments) were extensively discussed in \cite{winter2019geometry}.
We~recap~the relevant results and try to convey an intuition for how rigidity is connected to arrangement spaces and~represen\-tation theory.

\begin{remark}
Deformations of $\Sigma$-realizations as in \cref{def:rigid} preserve the symmetries in $\Sigma$.
A rigid realization cannot be deformed without loosing these relevant symmetries. 
However, such a rigid realization might well be \enquote{flexible} in the usual sense of rigidity theory of frameworks, where one cares mainly~about preserving edge-lengths. 
Preserving edge-lengths is also \mbox{not necessary~for~a~deforma}\-tion in our sense (unless, the symmetry requirements enforce it, see \cref{rem:edge_transitive_rigidity}).
%
\end{remark}

\begin{observation}\label{res:t_Ut}
Of course, any continuous reorientation or rescaling is a deformation in the sense of \cref{def:rigid}.
But more interesting are the deformations that acts in a \mbox{non-linear} way, as \eg\ the deformation of the skeleton of the truncated tetrahedron in \cref{ex:truncated_tetrahedron}.

Let us call such a non-linear deformation a \emph{proper deformation}.
We can formulate this in terms of arrangement spaces: a deformation $t\mapsto v(t)$ induces a map $t\mapsto U(t)\subseteq\RR^n$, that assigned to $t$ the arrangement spaces of $v(t)$.
Recalling \cref{res:equivalence}, a \enquote{proper deformation} is a deformation for which the map $t\mapsto U(t)$ is non-constant.
\end{observation}

\begin{observation}\label{res:deformation_intuition}
In a sense, the map $t\mapsto U(t)$ from \cref{res:t_Ut} describes a continuous curve in the \enquote{space of $d$-dimensional subspaces of $\RR^n$} (formally known as the \emph{Grassmannian} $G(d,n)$).
When transitioning continuously from a subspace $U(0)$ to a distinct subspace $U(1)$,
one necessarily passes though infinitely many pairwise distinct $\Sigma$-invariant subspaces $U(t)$.
This gives a necessary condition for the existence of a proper deformation: 

\begin{corollary}\label{res:finitely_many_subspaces}

If there are only finitely many $\Sigma$-irreducible subspaces, then all~$\Sigma$-realizations are rigid.
\end{corollary}

Note that the statement of \cref{res:finitely_many_subspaces} is slightly stronger than what we~have argued for in \cref{res:deformation_intuition}: in the corollary we require finitely many \emph{$\Sigma$-irreducible} subspaces rather than \emph{$\Sigma$-invariant} subspaces.
\begin{proof}[Proof of \cref{res:finitely_many_subspaces}]
%
Since every $\Sigma$-invariant subspace is the direct sum of $\Sigma$-irreducible subspaces, if there are only finitely many of the latter, then there are only finitely many of the former, and no $\Sigma$-realization can be flexible as discussed im \cref{res:deformation_intuition}.
\end{proof}

\end{observation}




In \cite{winter2019geometry} we have deduces several other necessary and sufficient conditions for~the~exis\-tence of deformations.
We list some of these:


\begin{theorem}\label{res:rigid}
\quad
\begin{myenumerate}
	\item If two irreducible realizations $v,w\in\mathcal R_d(G,\Sigma)$ have non-orthogonal arrangement spaces, then they can be deformed into each other, in particular, both realizations are flexible. (\!\cite{winter2019geometry}, Theorem 5.11)
	\item If an irreducible realization $v\in\mathcal R_d(G,\Sigma)$ is flexible, then $v$ can be deformed into a non-equivalent realization $w\in\mathcal R_d(G,\Sigma)$ with an arrangement space non-orthogonal to the one of $v$. (\!\cite{winter2019geometry}, Theorem 5.11)
	\item If two irreducible realizations $v,w\in\mathcal R_d(G,\Sigma)$ can be deformed into each other, then their representations are isomorphic, or if reoriented appropriately, we can assume that they have the same representation.\footnote{Together with point $(i)$, this is the realization version of the representations theoretic fact in \cref{res:non_orthogonal_subspaces}.}
	(\!\cite{winter2019geometry}, Corollary 5.5)
\end{myenumerate}
\end{theorem}

A first relevant application for us is the following:

\begin{theorem}\label{res:deform_into_balanced}
Every irreducible realization $v\in\mathcal R_d(G,\Sigma)$ can be deformed into a balanced $\Sigma$-realization.

Moreover, if $v$ is not already balanced, then $v$ can be deformed into at least two balanced $\Sigma$-realizations $w$ and $w'$ to different eigenvalues.
\begin{proof}
Consider the eigen-decomposition
$$\RR^n=\Eig_G(\theta_1)\oplus\cdots \oplus\Eig_G(\theta_m)$$
of $\RR^n$ into pairwise orthogonal $\Sigma$-invariant subspaces (the eigenspaces are $\Aut(G)$-invariant as seen the proof of \cref{res:spectral_then_symmetric}).

If $v$ were balanced, then we were done.
So assume, that $v$ is not balanced.
Its arrangement space $U\subseteq\RR^n$ is therefore not contained in any of the eigenspace of $G$.
Hence, $U$ is non-orthogonal to at least two of the eigenspaces, say $\Eig_G(\theta_1)$ and $\Eig_G(\theta_2)$.
Let $U_i:=\pi_i(U)$ denote the orthogonal projection of $U$ onto $\Eig_G(\theta_i)$.
Since $v$ is irreducible, so is $U$ as a subspace of $\RR^n$ (by \cref{res:arrangement_space_symmetric} $(ii)$), and one can check that $U_i$~is itself $\Sigma$-irreducible and non-orthogonal to $U$ (or see \cref{res:non_orthogonal_projection} in the appendix).
Any $\Sigma$-realization $w^{\smash{(i)}}$ with arrangement space $U_i\subseteq\Eig_G(\theta_i)$ is now irreducible (by \cref{res:arrangement_space_symmetric} $(ii)$), balanced (by \cref{res:arrangement_space_spectral}), and can be deformed into $v$ (by \cref{res:rigid} $(i)$).
\end{proof}
\end{theorem}


\begin{corollary}\label{res:rigid_implies_balanced}
If $v\in\mathcal R_d(G,\Sigma)$ is rigid and irreducible, then it is balanced.
\end{corollary}

\begin{theorem}\label{res:eigenvalue_multiplicity_implications}
Let $v\:V\to\RR^d$ be an irreducible $\Sigma$-realization of $G$. Suppose~that $G$ has a single eigenvalue $\theta$ of largest multiplicity $\mu_1$, and the second largest multiplicity of any eigenvalue of $G$ is $\mu_2$. Then
\begin{myenumerate}
	\item if $\mu_2<d$, then $v$ is balanced with eigenvalue $\theta$.
	\item if additionally $\mu_1 <2d$, then $v$ is rigid.
\end{myenumerate}
\begin{proof}
If $v$ is not already balanced, then by \cref{res:deform_into_balanced} it can be deformed into two balanced realizations with different eigenvalues.
Then these realizations~must~be of dimension $d$, and their arrangement spaces must be contained in different eigen\-spaces. But if $\mu_2<d$, then only one of these eigenspaces can have a large enough dimension.
Thus, $v$ must already have been balanced, and its arrangement space $U\subseteq\RR^n$ must be contained in the only large enough eigenspace $\Eig_G(\theta)$.
This proves $(i)$.

If $v$ were flexible, then it can be deformed into another non-equivalent irreducible $\Sigma$-realization $w$, also of dimension $d$.
For the same reasons as before, $w$ would be balanced with eigenvalue $\theta$, and thus its arrangement space $U'\subseteq\RR^n$ would be contained in $\Eig_G(\theta)$.
In particular, we have $U+U'\subseteq\Eig_G(\theta)$.
But since $U$ and $U'$ are irreducible, their intersection is trivial, and we have $$\mu_1=\dim\Eig_G(\theta)\ge \dim(U+U')=\dim(U)+\dim(U)-\dim(U\cap U')=2d.$$
Thus, if $\mu_1<2d$, then $v$ must have been rigid, proving $(ii)$.
\end{proof}
\end{theorem}


\subsection{Summary}
The take away messages of this section are:

\begin{itemize}
\item
A realization being symmetric, spectral, balanced, irredu\-cible or rigid can be nicely encoded in the language of arrangement spaces, which makes this concept an attractive tool for us.

\item Spectral realizations are as symmetric as possible. 
Since they are also fast to compute, they present a useful construction, even for those, mainly~interes\-ted in symmetric realizations.

\item
But, there are (irreducible) symmetric realizations that are not spectral,~only balanced, or not even that.
So in general, spectral methods are not enough.

\item Still, every irreducible realization can be deformed into a balanced realization (while keeping its symmetries). Conversely, this also means, that if a realization is rigid, then it must be balanced.
\end{itemize}

\tempnewpage

\section{Vertex-transitive realizations}
\label{sec:vertex_transitive}

Let us formally introduce what we mean when we say that a realization is vertex-transitive.

\begin{definition}\label{def:vertex_transitive}
Given a graph $G$.
\begin{myenumerate}
	\item $\Sigma\subseteq\Aut(G)$ is called \emph{vertex-transitive} if it acts transitively on $V$.
	\item $G$ is called \emph{vertex-transitive} if $\Aut(G)$ is vertex-transitive.
	\item A realization is called \emph{vertex-transitive} if it is a $\Sigma$-realization for some~vertex-tran\-sitive $\Sigma\subseteq\Aut(G)$.
\end{myenumerate}
\end{definition}

The notions in points $(ii)$ and $(iii)$ can be adopted for all other kinds of transitivities that we will encounter later on, as \eg\ edge- or arc-transitivity.

\begin{observation}\label{res:vertex_transitive}
In a vertex-transitive realization $v$, every vertex can be mapped onto any other vertex by an orthogonal transformation.
Thus, all vertices must be on a common sphere of radius $r(v)$ around the origin.
If the realization is norma\-lized, equation \eqref{eq:radius} states that for all $i\in V$ holds
\begin{equation}\label{eq:radius2}
\|v_i\|^2=[r(v)]^2=\frac dn.
\end{equation}
\end{observation}

\begin{observation}\label{res:single_vertex}
A vertex-transitive $\Sigma$-realization $v$ is already determined by its representation $T$, and the placement of a single vertex, say $v_1$.
Since for every $i\in V$ there is a $\sigma\in \Sigma$ with $\sigma(1)=i$, we have
$$v_i=v_{\sigma(1)}=T_\sigma v_1.$$
Note however, that not all choices of $v_1\in\RR^d$ are feasible (see also \cref{res:fix_rigid}, or Construction 4.4 in \cite{winter2019geometry}).
\end{observation}

The rest of this section is devoted to a discussion on how vertex-transitivity (or the absence of it) influences rigidity.
Recall, that we already met vertex-transitive realizations that are rigid (the dodecahedron in  \cref{ex:dodecahedron}, rigidity will be proven in \cref{sec:distance_transitive}), or that are flexible (the truncated tetrahedron in \cref{ex:truncated_tetrahedron}).


We start out with a rigidity-criterion specifically for vertex-transitive \mbox{realizations:}

\begin{lemma}
\label{res:fix_rigid}
Given a vertex-transitive realization~$v\in$ $\mathcal R_d(G,\Sigma)$ with representa\-tion \mbox{$T\:\Sigma\to\Ortho(\RR^d)$}. Let $\Sigma_i:=\{\sigma\in\Sigma\mid \sigma(i)=i\}$ denote the stabilizer of $\Sigma$~at~$i\in V$\!, and define
$$\Fix(T,\Sigma_i):=\{x\in\RR^d\mid T_\sigma x=x \text{ for all $\sigma\in\Sigma_i$}\}\subseteq\RR^d.
$$
%
If $\dim\Fix(T,\Sigma_i)=1$ for some (and then all) $i\in V$, then $v$ is rigid.
\begin{proof}
By $\dim\Fix(T,\Sigma_i)=1$ and $v_i\in\Fix(T,\Sigma_i)$ we have $\Fix(T,\Sigma_i)=\Span\{v_i\}$.

Suppose that $v$ can be deformed into $w\in\mathcal R_d(G,\Sigma)$.
Then (by an appropriate reorientation) both realization can be assumed to have the same \mbox{representation} $T$ (\cref{res:rigid} $(iii)$).
But then, for the same reason as for $v_i$, we have \mbox{$w_1\in\Fix(T,\Sigma_1)$} $=\Span\{v_1\}$, \ie\ $w_1=\alpha v_1$ for some $\alpha\in\RR$.
By vertex-transitivity (see \cref{res:single_vertex}) this holds for all vertices, that is, $v=\alpha w$. 

Hence $v$ cannot be deformed into a non-equivalent realization, and is rigid.
\end{proof}
\end{lemma}

We will make use of this in the proof of \cref{res:full_local_dimension_is_good}, when we study properties of \emph{arc-transitive} realizations.

From \cref{res:fix_rigid}, we immediately have that all 1-dimensional vertex-transitive realizations are rigid.
However, the converse of \cref{res:fix_rigid} is not true, that is, that there are rigid realizations with $\dim\Fix(T,\Sigma_i)$ $\ge 2$ for all $i\in V$.

\begin{example}\label{ex:exception}
This phenomenon is not primarily linked to graph realizations~(it already occurs for point arrangements), but, for better visualization, we can demonstrate its effect on the cycle graph $C_n$ on $n\ge 3$ vertices $V=\{1,...,n\}$.

Consider the cyclic subgroup $\Sigma\subset\Aut(C_n)$ generated by the one-cycle permuta\-tion \mbox{$\sigma=(12\cdots n)\in\Sigma$}.
The group $\Sigma$ acts freely on $V\!$, that is, $\Sigma_i=\{\id\}$ for~all~verti\-ces $i\in V$.

Let $T\:\Sigma\to\Ortho(\RR^2)$ be the $\Sigma$-representation that maps \mbox{$\sigma\mapsto R_{2\pi/n}\in\Ortho(\RR^2)$ -- the} rotation of the plane by $2\pi/n$ around the origin. By $\Sigma_i=\{\id\}$ we necessarily have $\Fix(T,\Sigma_i)=\RR^2$ for all vertices $i\in V$.

Nevertheless, all non-zero realizations with the representation $T$ give the skeleton of the regular $n$-gon in different orientations,
\begin{center}
\includegraphics[width=0.45\textwidth]{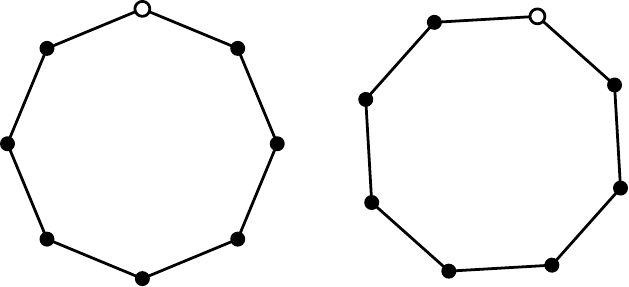}
\end{center}
and so all the realizations to this representation are equivalent, and no deformation is proper.
\end{example}

It would be interesting to determine the complete list of all exceptions, that is, of all rigid realizations with $\dim\Fix(T,\Sigma_i)\ge 2$ (see \cref{q:dim_fix_ge_2}).

We close this section with a note on what happens in the absence of vertex-transi\-tivity.
We show that for a not vertex-transitive graph all \enquote{interesting} \mbox{$\Sigma$-reali}\-zations (in a sense formalized in \cref{res:vertex_transitive_rigid}) are flexible.
This is because the vertices in each orbit can be placed independently of each other, creating a degree of freedom and preventing rigidity.

\begin{observation}
\label{res:vertex_transitive_rigid}
%
If the group $\Sigma$ does not act transitively on $V$, then $$V=V_1\cupdot\cdots\cupdot V_m$$ \mbox{decomposes} into $m\ge 2$ $\Sigma$-orbits $V_k$.
Now, if $v\: V\to\RR^d$ is a $\Sigma$-realization and $\alpha=(\alpha_1,...,\alpha_m)\in\RR^m$ is an $m$-tuple of real numbers, then the realization $v(\alpha)\: V\to\RR^d$, with~$v_i(\alpha):=$ $\alpha_k v_i$ whenever $i\in V_k$, is also~a~$\Sigma$-realization with the same repre\-sentation.
A continuous function $\alpha(t)\:[0,1]\to $ $\RR^m$ with $\alpha(0)=(1,...,1)$ induces a deformation $v(\alpha(t))$ from $v$ into $v(\alpha(1))$ (at least it does so under some conditions, that we discuss below).

For example, the deformation of the 6-cycle $C_6$, shown below, is of this form.~The group $\Sigma\subset\Aut(G)$ is choosen to have two orbits on $V$, indicated by the colors of the points. 
On the right, one of the orbits collapses to a single point in the origin.
\begin{center}
\includegraphics[width=0.8\textwidth]{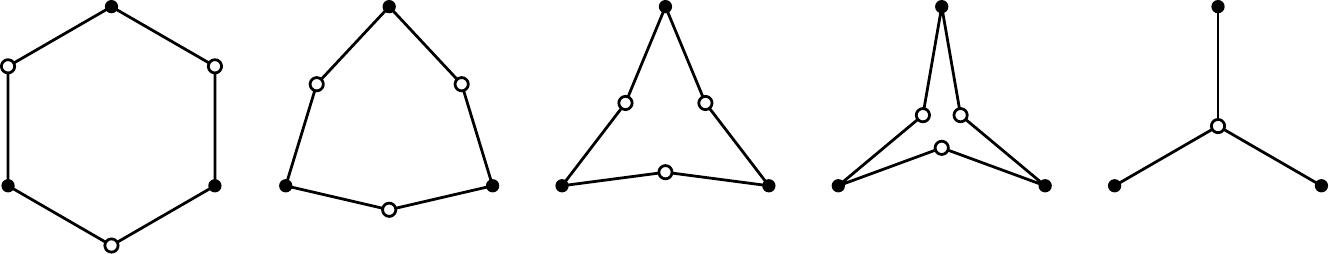}
\end{center}

Now, does this imply that a not vertex-transitive realization is always \mbox{flexible?~Not} necessarily: the deformation induced by $\alpha(t)$ might not be proper, as \eg\ in the~case of the \enquote{rhombus realization} of the 4-cycle $C_4$ shown below:
\begin{center}
\includegraphics[width=0.65\textwidth]{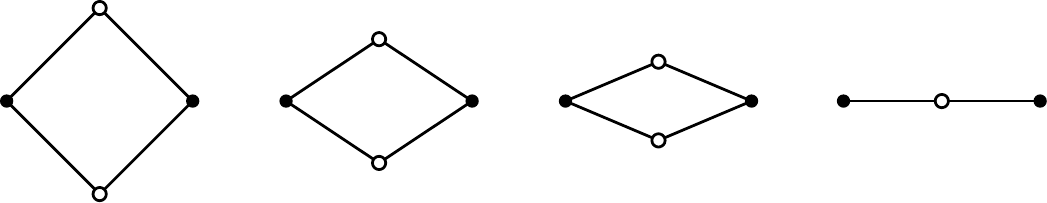}
\end{center}
Again, this is a deformation of the form discussed above, but this time, all $v(\alpha(t))$ are equivalent.
Note further, that the realization on the right is not part of this deformation as it is not full-dimen\-sional (recall that $\mathcal R_d(G,\Sigma)$ contains only full-dimensional realizations).
In~fact,~these $\Sigma$-realizations of $C_4$ are reducible, and~the right-most image shows one of the irreducible components of this $\Sigma$-realization instead.

Let us try to understand what distinguishes these two examples.
Suppose that $v$ is irreducible with representation $T$.
Define the sets $\mathcal N_k:=\{v_i\mid i\in V_k\}$ for all $k\in\{1,...,m\}$ (the image of the orbit $V_k$ under the realization $v$).
If $\mathcal N_k =\{0\}$ (as for the \enquote{white} orbits in the two right images), then~this is equivalent to ignoring the orbit $V_k$, and we could have studies $V\setminus  V_k$ instead.
Let us therefore assume that all these sets are non-zero.
Note that $\Span (\mathcal N_k)$ is a $T$-invariant subspace of $\RR^d$.
Since $v$ is irreducible and $\mathcal N_k$ is non-zero, we necessarily have $\Span(\mathcal N_k)=\RR^d$.
But then, if a transformation of $\RR^d$ fixes one of the $\mathcal N_k$ point-wise, this transformation must be the identity, hence fixes all $\mathcal N_k$ point-wise.
For such realizations, deformations of the form $v(\alpha(t))$ are indeed proper, and $v$ is flexible.
\end{observation}

\begin{corollary}
\label{res:vertex_transitive_rigid_cor}
If an irreducible realization is not vertex-transitive but has at least two non-zero orbits on the vertices (the $\mathcal N_i$ in \cref{res:vertex_transitive_rigid}), then the realization is flexible.
\end{corollary}

The quintessence of \cref{res:vertex_transitive_rigid} and \cref{res:vertex_transitive_rigid_cor} is meant to be the following: the study of rigidity is much more interesting for vertex-transitive realizations, and we shall therefore focus on these.

\tempnewpage

\section{Edge- and arc-transitive realizations}
\label{sec:arc_transitive}

In this section we explore the properties of edge- and arc-transitive realizations.
Recall that an \emph{arc} in a graph is an incident vertex-edge pair.
The notions of edge- and arc-transitive graphs and realizations are defined parallel to \cref{def:vertex_transitive}.

\begin{observation}\label{res:edge_transitive}
If $v$ is edge-transitive, then all edges can be mapped onto each other by orthogonal transformations.
Thus, all edges have the same length,~and~their end vertices have the same inner product.
That is, the following notions are well-defined:
$$\omega(v):=\<v_i,v_j\>,\qquad\ell(v):=\|v_i-v_j\|$$
for any $ij\in E$.
The latter is called \emph{edge length} of $v$.
\end{observation}

\begin{remark}\label{rem:edge_transitive_rigidity}
Let $v(t)\:[0,1]\to\mathcal R_d(G,\Sigma)$ be a deformation for some edge-transitive $\Sigma\subseteq\Aut(G)$.
Then, the deformation $w(t):=v(t)/\ell(v(t))$ has edge length \mbox{$\ell=1$~for~all} $t\in[0,1]$.

In this sense, a deformation \wrt\ an edge-transitive groups can always be consi\-dered as a flex in the usual sense of rigidity theory, which preserves edge-lengths.
\end{remark}

Under certain conditions, the quantities from \cref{res:edge_transitive} can be computed~ex\-plicitly.

\begin{proposition}\label{res:cosine_and_length}
If $v$ is normalized, vertex- and edge-transitive, and $\theta$-balanced~for some eigenvalue $\theta\in\Spec(G)$ (\resp\ Laplacian eigenvalue $\lambda :=\deg(G)-\theta$, see also \cref{sec:appendix_spec_theory}), then
\begin{equation}\label{eq:cosine_and_length}
\omega(v)=\frac{\theta d}{2|E|},\qquad [\ell(v)]^2=\frac{\lambda d}{|E|}.
\end{equation}
\begin{proof}
Using equations \eqref{eq:balanced} and \eqref{eq:radius2} we see that for all $ij\in E$ holds
$$
\deg(G)\cdot\omega(v)
=\sum_{\mathclap{j\in N(i)}} \<v_i,v_j\>
=\Big\< \!v_i,\sum_{\mathclap{j\in N(i)}} v_j\Big\> \overset{\text{\eqref{eq:balanced}}}= \<v_i,\theta v_i\> 
=\theta \cdot \|v_i\|^2
\overset{\text{\eqref{eq:radius2}}}= \frac{\theta d}n.
$$
Expression \eqref{eq:cosine_and_length} for $\omega(v)$ then follows by $\deg(G)\cdot n=2\card E$.
The expression for $\ell(v)$ can now be derived as follows:
%
\begin{align*}
[\ell(v)]^2 
&= \|v_i-v_j\|^2 
= \|v_i\|^2+\|v_j\|^2-2\<v_i,v_j\> 
\\&= 2\big( [r(v)]^2-\omega(v)\big)
= 2\Big(\frac dn-\frac{d\theta}{2|E|}\Big) = 2d\Big(\frac{\deg(G)}{2|E|}-\frac\theta{2|E|}\Big)
\\&= \frac{2d(\deg(G)-\theta)}{2|E|} = \frac{\lambda d}{|E|}.
\end{align*}
%
%
\end{proof}
\end{proposition}

If $v$ is spherical instead of normalized, one can still compute the follow\-ing relative 
quantities
\begin{equation}\label{eq:relative}
\frac{\omega(v)}{[r(v)]^2} = \frac{\theta}{\deg(G)},\wideand \Big[\frac{\ell(v)}{r(v)}\Big]^2 = \frac{2\lambda}{\deg(G)} = 2\Big(1-\frac\theta{\deg(G)}\Big),
\end{equation}
named \emph{cosine} and \emph{relative length} of $v$.
The name \enquote{cosine} is because for $ij\in E$~holds
\begin{equation}\label{eq:angle}
\cos \angle(v_i,v_j) = \frac{\<v_i,v_j\>}{\|v_i\|\|v_j\|} = \frac{\omega(v)}{[r(v)]^2} = \frac{\theta}{\deg(G)}.
\end{equation}

We can use these formulas to perform some quick computations on some special poly\-topes.
These are quite laborious if done by hand.

\begin{example}\label{ex:circumradius}
In \cref{ex:dodecahedron}, we have seen that the skeleton of the dodecahedron can be obtained as the $\theta_2$-realization of its edge graph (where $\theta_2=\sqrt 5$).
We can then use \eqref{eq:relative} to compute the circumradius of the dodecahedron with edge length $\ell=1$:
$$r(v)= \sqrt{\frac{\deg(G)}{2\lambda_2}} = \sqrt{\frac{3}{2(3-\sqrt{5})}}\approx 1.401258.$$
\end{example}

\begin{example}\label{ex:dihedarl_angle}
The \emph{dihedral angle} of a $d$-dimensional polytope $P\subset\RR^d$ at a $(d-2)$-dimensional face $\sigma$ (also called a \emph{ridge} of $P$) is the angle between the two facets incident to $\sigma$, when measured on the inside of the polytope.
In the 3-dimensional case, the ridges are just the edges of $P$.

For example, for the icosahedron (the dual of the dodecahedron), this angle is the same for each edge.
The dihedral angle $\alpha$ of a the icosahedron is exactly $\pi$~minus the angle between incident vertices in its dual. 
The angle between vertices can be computed via \eqref{eq:angle}, and we find:
$$\alpha = \pi -\arccos\Big( \frac{\theta_2}{\deg(G)}\Big) = \pi-\arccos\Big(\frac{\sqrt 5}{3}\Big)\;\hat\approx\; 138.1896^\circ.$$

\begin{center}
\includegraphics[width=0.5\textwidth]{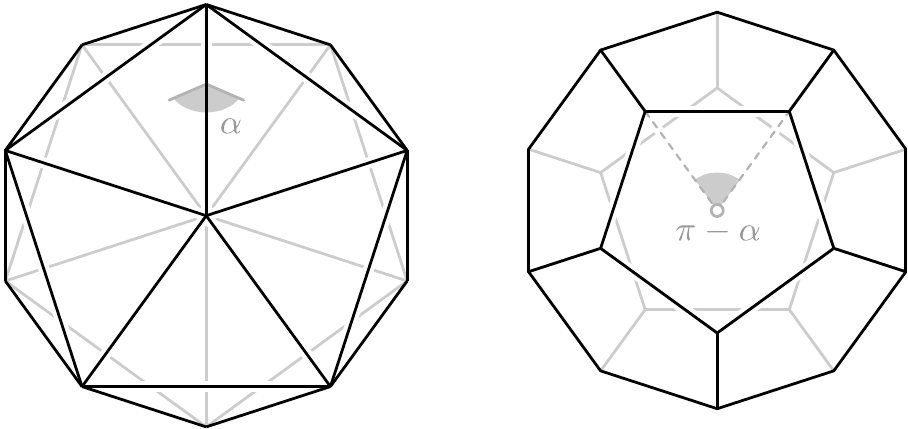}
\end{center}
\end{example}

The computations in \cref{ex:circumradius} and \cref{ex:dihedarl_angle} work equivalently for all the other regular polytopes as noted in the last paragraph of \cref{ex:dodecahedron}.



We have seen that general symmetric realizations are not necessarily spectral, balanced or rigid.
This includes vertex-transitive realizations as we have seen~in~\cref{ex:truncated_tetrahedron}.
We also discussed that vertex-transitivity appears as a plausible minimal require\-ment for obtaining rigidity result (see \cref{res:vertex_transitive_rigid}).
In particular, purely edge-transitive realizations are not expected to be either rigid or balanced.

\begin{example}
The skeleton of the (edge-transitive) \emph{rhombic dodecahedron} (a \emph{Catalan solid}, see the middle image below), is an $\Aut(G)$-realization of its edge-graph.
The image below shown a deformation of this realization (of the form, as constructed in \cref{res:vertex_transitive_rigid}):
\begin{center}
\includegraphics[width=0.65\textwidth]{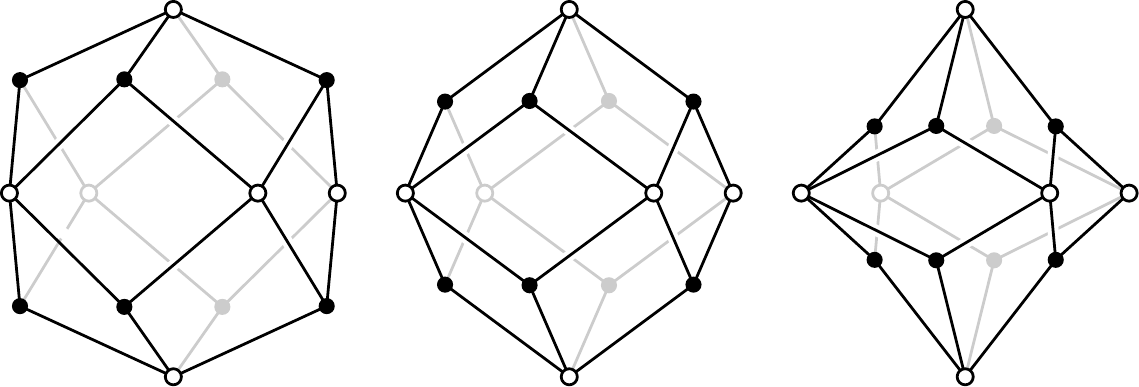}
\end{center}
It was also realized by Licata and Powers \cite{licata1986surprising}, that the $\theta_2$-realization of the edge-graph of the rhombic dodecahedron gives exactly its skeleton, rather than any of the infinitely many other deformations.
This is because its skeleton is balanced,~but in contrast to the regular polytopes, this cannot be explained from general symmetry considerations (see \cref{res:full_local_dimension_is_good} and \cref{ex:24_cell}) and appears more accidental.
\end{example}

For the reasons presented before, we from now on also assume vertex-transitivity.
While there is a difference between being simultaneously vertex- and edge-transitive and being arc-transitive (the latter is strictly stronger, in between these classes there exist the so-called \emph{half-transitive graphs}, see \cite{bouwer1970vertex, holt1981graph}), in the following we primarily focus on arc-transitive graphs and realizations.

As all the symmetry classes before, also arc-transitive realizations are not necessarily rigid or spectral, though it becomes increasingly more complicated to construct counterexamples.

\begin{example}\label{ex:K_nn}
Consider the complete bipartite graph $K_{n,n}$ on $2n$ vertices,~which is arc-transitive. 
Its spectrum is $\{(-n)^1,0^{2(n-1)},n^1\}$, and indeed, the 0-eigen\-space~is $\Aut(K_{n,n})$-irreducible.
But there are many arc-transitive subgroups of $\Aut(K_{n,n})$ for which this~eigen\-space decomposes into smaller irreducible subspaces, giving rise to non-spectral arc-transitive realization of $K_{n,n}$.
Many of these are also flexible.
We describe a general procedure to construct such.

Let $v(t)\:[0,1]\to\mathcal R_d(G,\Sigma)$ be proper deformation between irreducible vertex-transitive $\Sigma$-realization of some graph $G$ on the vertex set $V=\{1,...,n\}$ (\eg\ the graph from \cref{ex:truncated_tetrahedron}).
Consider $K_{n,n}$ with vertex set $V_1\cupdot V_2$, where the $V_i$ are disjoint copies of $V$.
Let $\Sigma_i\subseteq \Aut(K_{n,n})$ be an isomorphic copy of $\Sigma$ acting on $V_i$ instead of $V$, and let $\tau\in\Aut(K_{n,n})$ be the involution that exchanges the both partition classes in the obvious way.
Then set $\Sigma':=\<\Sigma_1,\Sigma_2,\tau\>\subseteq\Aut(K_{n,n})$.
Note first, that $\Sigma'$ indeed acts arc-transitively on $K_{n,n}$.
We further claim that the 0-eigenspace of $K_{n,n}$ is \emph{$\Sigma'$-reducible}, and that some irreducible subspace corresponds to a \emph{flexible} arc-transitive realizations.

To see this, consider the following deformation $w(t)\: [0,1]\to\mathcal R_{2d}(K_{n,n},\Sigma')$:
$$w_i(t) = \begin{cases} (v_i(t),0) & \text{if $i\in V_1$}\\ (0,v_i(t)) & \text{if $i\in V_2$}\end{cases}.$$
Since $v(\free)$ is a proper deformation, so is $w(\free)$.
One checks that for all $t\in[0,1]$ $w(t)$ is indeed an irreducible $\Sigma'$-realization, which, by the existence of the deformation, must be flexible.

Since $d\ge 2$ and $K_{n,n}$ has only a single non-simple eigenvalue, from \cref{res:eigenvalue_multiplicity_implications} $(i)$ follows that $w(t)$ must be balanced with eigenvalue $0$. 
Since $w(t)$ is flexible, the premise of \cref{res:eigenvalue_multiplicity_implications} $(ii)$ cannot be satisfied.
The multiplicity of the eigenvalue $0$ must therefore be at least $4d$ (twice the dimension of $w(t)$), and so $w(t)$ cannot be spectral, but must only balanced.
\end{example}

The trick in \cref{ex:K_nn} was to choose $\Sigma\subset\Aut(K_{n,n})$ as not the full symmetry group of $K_{n,n}$, but as a proper subgroup.
It is unclear whether this was necessary, that is, whether an $\Aut(G)$-realization of an arc-transitive graph is always rigid (see \cref{q:arc_transitive_balanced}).
An $\Aut(G)$-realization of an arc-transitive graph needs not be spectral though, as we shall see in \cref{ex:C6_C6} further below.
All realizations in \cref{ex:K_nn} are balanced, and it is equally unclear, whether this is the case for all irreducible arc-transitive realizations (see \cref{q:arc_transitive_rigid}).

In the next section, we restrict to a class of arc-transitive realizations for which these questions can be resolved.

\subsection{Full local dimension}
\label{sec:full_local_dimension}

In this section we focus on  realizations of \emph{full local dimension}, which means that the edge directions at each vertex span the whole space.
For example, this is always the case for the skeleta of convex polytopes.

\begin{definition}
\label{def:full_local_dimension}
A realization $v$ is said to be of \emph{full local dimension} if 
$$\rank\{v_j-v_i\mid j\in N(i)\}=d,\quad\text{for all $i\in V$}.$$
\end{definition}

Being of full local dimensional implies being of full dimension,
but not every~full-dimensional realization is of full local dimension, not even if it is arc-transitive.

\begin{example}\label{ex:not_full_local_dim}
Consider the 4- and 5-dimensional spectral realizations of the edge-graph of the dodecahedron (we have seen in \cref{ex:dodecahedron} that this graph indeed has eigenvalues of multiplicity four and five).
As spectral realizations, they are~arc-transitive.
We will later see (in \cref{res:irreducibe_distance_transitive_is_spectral}) that these realizations are even~irre\-ducible.
However, they cannot be of full \emph{local} dimension since the edge-graph of the dodecahedron is only of degree three.
\end{example}

In general, a realization of full local dimension must neither be rigid/balanced (see \cref{ex:truncated_tetrahedron}) nor irreducible (\eg\ the skeleton of a prism).
This is different in the case of arc-transitive realizations.


\begin{theorem}\label{res:full_local_dimension_is_good}
Let $v$ be an arc-transitive realization of full local dimension. Then
\begin{myenumerate}
	\item $v$ is irreducible,
	\item $v$ is rigid, and
	\item $v$ is balanced.
\end{myenumerate}
%
%
\begin{proof}
Suppose $v$ is a $\Sigma$-realization for some arc-transitive $\Sigma\subseteq\Aut(G)$ and $T\:\Sigma\to\Ortho(\RR^d)$ its representation.

Let $\Sigma_i\subseteq\Sigma$ denote the stabilizer of $\Sigma$ at the vertex $i\in V$.
We want to determine the invariant subspaces of the restriction $T|_{\Sigma_i}\:\Sigma_i\to\Ortho(\RR^d)$.
Clearly, $\Span\{v_i\}$ is invariant, and $T|_{\Sigma_i}$ acts on it by identity (we say, it acts trivially).
In other words, $\Span\{v_i\}\subseteq\Fix(T,\Sigma_i)$ (as defined in \cref{res:fix_rigid}).
We show that we actually have equality.

By arc-transitivity, $T|_{\Sigma_i}$ acts transitively on the set $\mathcal N_i:=\{v_j\mid j\in N(i)\}$.
That is, for any two $w_1,w_2\in \mathcal N_i$ there exists a $\sigma\in\Sigma_i$ with $T_\sigma w_1=w_2$.
And so for any $x\in\Fix(T,\Sigma_i)$ (and by using that $T_\sigma$ is orthogonal) we have
$$\<x,w_1\> = \<T_\sigma x,T_\sigma w_1\> = \<x,w_2\> \quad\implies\quad \<x,w_1-w_2\> = 0.$$
And since this holds for all $x\in\Fix(T,\Sigma_i)$ and all pairs $w_1,w_2\in\mathcal N_i$, we obtained $\Fix(T,\Sigma_1) \subseteq \aff (\mathcal N_i)^\bot$.~So we found $\Span\{v_i\}\subseteq\Fix(T,\Sigma_i)\subseteq\aff(\mathcal N_i)^\bot$.
But from full local \mbox{dimension follows} $\dim\aff(\mathcal N_i)\ge d-1\implies \dim\aff(\mathcal N_i)^\bot\le 1$.
Thus, the dimensions of the~sub\-spaces in the inclusion chain must agree, and we actually have $\Span\{v_i\}=\Fix(T,\Sigma_i)=\aff(\mathcal N_i)^\bot$.
In particular, we have $\dim\Fix(T,\Sigma_i)=1$ for~all $i\in V$, and so $v$ is rigid by \cref{res:fix_rigid}. This proves $(ii)$.


To proceed, we show that all $T|_{\Sigma_i}$-invariant subspace, besides $\Span\{v_i\}$, are~contained in the orthogonal complement $v_i^\bot$.
So suppose that $W\not\subseteq v_i^\bot$ is a \mbox{$T|_{\Sigma_i}$-irre}\-ducible subspace, not contained in the orthogonal complement of $v_i$.
In other words, $W$ is non-orthogonal to $\Span\{v_i\}$. 
But if irreducible subspaces are non-orthogonal, then the representation acts isomorphically on them (see \cref{res:non_orthogonal_subspaces}).
So, since $T|_{\Sigma_i}$ acts trivially on $\Span\{v_i\}$, it must act trivially on $W$ too.
But then $W\subseteq\Fix(T,\Sigma_i)=\Span\{v_i\}$.

We can now show that $v$ is irreducible.
For this, suppose that $W\subseteq\RR^d$ is a~$T$-invariant subspace of $\RR^d$.
Such a subspace must also be invariant \wrt\ all the $T|_{\Sigma_i}$, that is, for each $i\in V$ it must either contain $\Span\{v_i\}$, or must be contained in $v_i^\bot$.
Because of vertex-transitivity, if $\Span\{v_i\}\subseteq W$ for some $i\in V$, then for all $i\in V$. 
Since $v$ is full-dimensional, this would give $W=\RR^d$.
Likewise, if $W\subseteq v_i^\bot$ for some $i\in V$, then for all $i\in V$.
Since $v$ is full-dimensional, this would give $W=\{0\}$.
Thus, $W$ is a trivial invariant subspace, and $v$ is irreducible.
This proves $(i)$.

Finally, since $v$ is rigid and irreducible, it follows that $v$ is balanced by \cref{res:rigid_implies_balanced}, which proves $(iii)$.
\end{proof}
\end{theorem}


\begin{corollary}
If an arc-transitive realization is reducible, then it cannot be of full local dimension.
\end{corollary}

In particular, these results apply to the skeleta of arc-transitive polytopes.
We discuss this in the case of regular polytopes:


\begin{example}\label{ex:24_cell}
\Cref{res:full_local_dimension_is_good} now assures us that the skeleta of regular polytopes (actually, of arc-transitive polytopes) are balanced. 
As mentioned in \cref{ex:dodecahedron}, they are actually spectral, which was shown in \cite{licata1986surprising} for all regular polytopes excluding the exceptional 4-dimensional regular polytopes.
We fill in this gap now.

Since we can now apply  \eqref{eq:relative} (which requires being balanced), it remains purely a matter of checking tabulated values for circumradii and edge-lengths of regular polytopes, as well as eigenvalues of their edge-graphs, to find that their skeleta are indeed $\theta_2$-realizations.

We demonstrate this on the example of the 24-cell, whose edge-graph we shall~call $G$.
One of the notable properties of the 24-cell is that the edge length equals~the~cir\-cumradius, or $\ell(v)/r(v)=1$ for its skeleton $v$.
Since $v$ is an arc-transitive realization, \ie~balanced by \cref{res:full_local_dimension_is_good}, we know that \eqref{eq:relative} applies.
We \mbox{can~rearrange} \eqref{eq:relative} for $\theta$ to find
\begin{equation}\label{eq:theta_from_metrics}
\theta=\deg(G)\Big(1-\frac12\cdot\Big[\frac{\ell(v)}{r(v)}\Big]^2\Big)=8\cdot\Big(1-\frac12\Big)=4.
\end{equation}
Indeed, the spectrum of $G$ is $\{(-8)^1,(-4)^4,0^6,4^4,8^1\}$ with eigenvalue $\theta_2=4$~of~multi\-plicity \emph{four}.
So the skeleton of the 24-cell is spectral with eigenvalue $\theta_2$.

The same technique works verbatim for all other regular polytopes as well.
The table below lists the degree, circumradius and edge-length of the three exceptional 4-dimensional regular polytopes, as well as the second largest eigenvalue $\theta_2$ for their edge-graphs (which have multiplicity four in all cases).
The reader can then check~that these values satisfy an equation analogue to \eqref{eq:theta_from_metrics}, \ie\ the equation produces the second-largest eigenvalue.
\begin{center}
\begin{tabular}{l|l|l|l|l}
polytope & $\deg(G)$ & $\ell(v)$ & $r(v)$ & $\theta_2$ \\
\hline
24-cell & \phantom08 & 1 & 1 & 4 \\
120-cell & \phantom04 & $3-\sqrt 5$ & $\sqrt 8$ & $2\phi-1$ \\
600-cell & 12 & $\phi^{-1}$ & $1$ & $2(1+\sqrt 5)$
\end{tabular}
\end{center}
Here, $\phi:=(\sqrt 5+1)/2$ denotes the \emph{golden ratio}.
Despite that we can check this on~a case-by-case basis, it remains a mystery why all these skeleta are spectral realizations.
In contrast, for all regular polytopes excluding these 4-dimensional exception, a satisfying answer will be given in the next section, in \cref{ex:regular_licata}.

We shall give a second proof for the 24-cell in \cref{ex:24_cell_cuboctahedron} where we do not~need to know the spectrum of the edge-graph.
\end{example}

%
%

We close this section with an example of an arc-transitive realization of full-local dimension that is \emph{not} spectral, showing that \cref{res:full_local_dimension_is_good} $(iii)$ cannot be improves in general.

\begin{example}\label{ex:C6_C6}
Consider the graph $G:=C_6\times C_6$, best visualized as the edge-graph of the hexagonal torus:
\begin{center}
\includegraphics[width=0.4\textwidth]{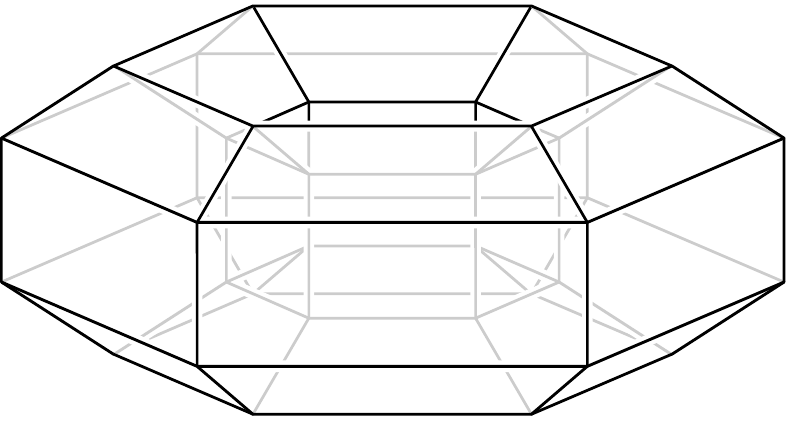}
\end{center}
While the torus is not arc-transitive, the graph $G$ is.
Its spectrum is 
$$\{(-4)^1,(-3)^4,(-2)^4,(-1)^4,0^{10},1^4,2^4,3^4,4^1\},$$
in particular, $G$ has an eigenvalue 0 of multiplicity ten.
The following formula~des\-cribes a 2-dimensional balanced $\Aut(G)$-realization with eigenvalue zero:
$$V(C_6)\times V(C_6)\ni(i,j)\;\mapsto\; \begin{pmatrix}
(-1)^i\\
(-1)^j
\end{pmatrix}\in\RR^2.$$
This realization is of full local dimension, but it is not spectral, since its dimension is smaller than 10.
\end{example}

%

\subsection{Summary}

The take away messages of this section are:

\begin{itemize}
\item If an arc-transitive realization is balanced, many metric properties can already be compute from only knowing $G$ and the eigenvalue.
\item Arc-transitive realizations seem to be better behaved than more general~(\eg\ vertex-transitive) realizations, but they still do not have to be rigid or spectral.
\item We do not know whether they have to be balanced.
\item If we assume full local dimension, everything becomes better, and we can now show that the realization is rigid, irreducible and balanced (but not necessarily spectral).
\end{itemize}

\tempnewpage

\section{Distance-transitive realizations}
\label{sec:distance_transitive}

In this final section of the paper we discuss a class of realizations with a particularly high degree of symmetry, called \emph{distance-transitivity}.

Distance-transitive graphs, and their generalizations, the \emph{distance-regular graphs}, 
form a class of graphs especially accessible by spectral graph theory methods.
The standard literature for these is the monograph by Brouwer, Cohen and Neumaier \cite{brouwer1989distance}.
The generic distance-regular graph has a trivial automorphism group, and thus, those are less relevant to this discussion, and we shall restrict to the \mbox{distance-tran}\-sitive graphs.
Nevertheless, the reader familiar with distance-regular graphs will recognize many of our arguments.

Throughout this section, let $\dist(i,j)$ denote the \emph{distance} between any two vertices $i,j\in V$,~\ie\ the length of the shortest path between $i$ and $j$.
The \emph{diameter} $$\diam(G):=\max_{i,j\in V}\dist(i,j)$$ of $G$ is the maximal distance between any two of its vertices.


\begin{definition}\label{def:distance_transitive}
A group $\Sigma\subseteq\Aut(G)$ acts \emph{distance-transitively} on $G$ if it acts~transitively on each of the sets
$$D_\delta:=\{(i,j)\in V\times V \mid \dist(i,j)=\delta\} ,\quad\text{for each $\delta\in\{0,...,\diam(G)\}$}.$$
%
\end{definition}

Distance-transitive graphs and realizations are defined parallel to \cref{def:vertex_transitive}.
Note that being arc-transitive is equivalent to being transitive on the set $D_1$.
Hence, distance-transitivity implies arc-transitivity.

\begin{example}
Complete graphs and cycle graphs are distance-transitive.
More~ge\-ne\-rally, the edge-graphs of the regular polytopes are distance-transitive, with the usual 4-dimensional exceptions. 
Even stronger, the skeleton of any of these regular polytopes is a distance-transitive realization of the edge-graph (this is easy to check for the simplices and $n$-cubes, and the only two remaining cases are the dodecahedron and the icosahedron, which can be checked by hand).
%

Other examples of distance-transitive graphs that are not necessarily edge-graphs are the Petersen-graph and all complete $r$-partite graphs $K_{n,...,n}$.
\end{example}

A very relevant observation concerning distance-transitive symmetry was already made in \cite{brouwer1989distance}:


\begin{theorem}[\!\cite{brouwer1989distance}, p.\ 137,  Proposition 4.1.11]
\label{res:eigenspaces_eq_irreducible}
For every distance-transitive group $\Sigma\subseteq\Aut(G)$, the $\Sigma$-irreducible subspaces of $\RR^n$ are exactly the eigenspaces of $G$.
\end{theorem}

We highlight again the stark contrast to the arc-transitive case, where eigenspaces are not necessarily irreducible (we were not even able to prove that every irreducible subspace is contained in an eigenspace).
The observation in \cref{res:eigenspaces_eq_irreducible} was also made specifically for the Petersen graph by Du and Fan in \cite{du2013graph}.

We shall give a partial proof for \cref{res:eigenspaces_eq_irreducible} in \cref{sec:cosine_vectors}, that is, we will~show that the eigenspaces are irreducible \wrt\ distance-transitive symmetry (see \cref{res:cosine_vectors_are_determined}).
The technique used there (namely, \emph{cosine vectors}) admits a generalization to not necessarily distance-transitive graphs.


\Cref{res:eigenspaces_eq_irreducible} in the form of realizations reads as follows:

\begin{theorem}\label{res:irreducibe_distance_transitive_is_spectral}
The spectral realizations of a distance-transitive graphs $G$ are exactly the irreducible distance-transitive realizations of $G$.
\end{theorem}

In other words, all symmetric realizations of a distance-transitive graph can be obtained by spectral methods.
We list some further consequences.

\begin{corollary}
\label{res:distance_transitive_consequences}
Given a distance-transitive realization $v$, then
\begin{myenumerate}
\item $v$ is rigid,
\item $v$ is an $\Aut(G)$-realization, and
\item the following are equivalent: $v$ being balanced, spectral and irreducible.
\end{myenumerate}
\begin{proof}
Suppose that $v$ is a $\Sigma$-realization with distance-transitive group $\Sigma\subseteq\Aut(G)$
and arrangement space $U\subseteq\RR^n$.

Since $G$ has only finitely many eigenspaces, by \cref{res:eigenspaces_eq_irreducible} there are only finitely many $\Sigma$-irreducible subspaces of $\RR^n$. By \cref{res:finitely_many_subspaces}, all $\Sigma$-realizations are rigid. This proves $(i)$.

For $(ii)$, first assume that $v$ is irreducible.
Then $v$ is spectral by \cref{res:irreducibe_distance_transitive_is_spectral}, and~therefore an $\Aut(G)$-realization by \cref{res:spectral_then_symmetric}.
%
%
Now, if $v$ is not irreducible, then we can consider a decomposition of $v$ into~irre\-ducible $\Sigma$-realizations $v^{(1)},...,v^{(k)}$ as in \cref{rem:irreducible_decomposition} (we assume that $v$ is appropriately oriented). 
Each of these irreducible constituents is now an $\Aut(G)$-realization of $G$ with an irreducible $\Aut(G)$-representation $$T^{(k)}\:\Aut(G)\to\Ortho(\RR^{d_k}).$$
From this we construct an $\Aut(G)$-representation $T\:\Aut(G)$ $\to\Ortho(\RR^d)$ for $v$ as in \eqref{eq:representation_decomposition}. Thus, $v$ is an $\Aut(G)$-realization, proving $(ii)$.

Finally, we prove $(iii)$.
Being spectral and being irreducible are equivalent by \cref{res:irreducibe_distance_transitive_is_spectral}. 
Also spectral implies balanced, and so we only have to prove the converse: if $v$ is balanced, then $U\subseteq\Eig_G(\theta)$ for some eigenvalue $\theta\in \Spec(G)$. Now $U$ contains at least one irreducible subspace, but $\Eig_G(\theta)$ already is irreducible. Thus $U=\Eig_G(\theta)$, $v$ is spectral, and we proved $(iii)$.
%
%
%
%
%
%
%
%
\end{proof}
\end{corollary}

%
%
%

Note especially part $(ii)$ which can be informally stated as follows:
the distance-transitive symmetries of $G$ cannot be geometrically separated from the other symmetries of $G$.
As soon as one tries to realize the distance-transitive symmetry of $G$, one automatically realizes all symmetries of $G$.

This is remarkable, as we have already encountered vertex-, edge- and even arc-transitive realizations that are not $\Aut(G)$-realizations, that is, they still leave some symmetries unrealized.
It appears as if distance-transitivity lies beyond a threshold, from which on symmetries can no longer be distinguished geometrically.

\begin{remark}\label{ex:regular_licata}
Licata and Powers \cite{licata1986surprising} have proven their result about regular polytopes (discussed in \cref{ex:dodecahedron} and \cref{ex:24_cell}) by explicit computation.

We can now obtain their result as follows: the skeleton of a regular polytope (not a 4-dimensional exception) is distance-transitive (in particular, arc-transitive) and of full local dimension, thus irreducible by \cref{res:full_local_dimension_is_good}. But if it is irreducible and distance-transitive, then~it must be spectral by \cref{res:irreducibe_distance_transitive_is_spectral}.
\end{remark}

\subsection{Cosine vector and cosine sequence}
\label{sec:cosine_vectors}

The remainder of this section is devoted to an idea that can be used to prove a part of \cref{res:irreducibe_distance_transitive_is_spectral}, but which also admits some interesting generalizations.

\begin{definition}
For a vertex-transitive realization $v$ and some vertex $i\in V$, the \emph{cosine vector} $u\in\RR^n$ of $v$ is the vector with components $u_j:=\<v_i,v_j\>$ for all $j\in V$.
\end{definition}

Vertex-transitivity ensures that this definition of the cosine vector is independent of our choice of $i\in V$ up to some coordinate permutation. We can therefore~assume that $u_i:=\<v_1,v_i\>$.
Recall further, that the cosine vector is contained in the arrangement space of $v$ as seen in  \cref{res:eigenvectors_from_realization}.

The central idea concerning the consine vector is explained in the following observation:

\begin{observation}\label{res:cosine_vectors}
Suppose that we are given a $\theta$-balanced $\Sigma$-realization $v$ and we want to know whether $v$ is spectral.
This is easy if we know the multiplicity of $\theta$, so suppose that we do not.

Let $U\subseteq\RR^n$ be the arrangement space of $v$.
If $v$ were not spectral, then we find that the subspace $\bar U:=U^\bot\cap\Eig_G(\theta)$ is non-zero and $\Sigma$-invariant.
The corresponding realization $\bar v$ with arrangement spaces $\bar U$ is then also a $\theta$-balanced $\Sigma$-realization.
If $u,\bar u\in\RR^n$ are the consine vectors of $v$ and $\bar v$ respectively, then $\<u,\bar u\>=0$ because they are contained in the orthogonal subspaces $U$ and $\bar U$.


The idea is to show that, in the right setting, being balanced $\Sigma$-realizations to the same eigenvalue is already so restrictive, that the corresponding cosine vectors have no chance to be orthogonal.
If this is the case, then we found that $\Eig_G(\theta)$ is $\Sigma$-irreducible and $v$ is the $\theta$-realization.
\end{observation}

This idea most directly applies to distance-transitive realizations.
For these, one shows that the cosine vector depends only on $G$ and $\theta$ (\cref{res:cosine_vectors_are_determined}). 

\begin{observation}
If $v$ is distance-transitive, then the value of $u_i=\<v_1,v_i\>$ depends only on $\delta:=\dist(1,i)$.
One therefore groups all entries with the same distance to $1\in V$ and writes $u_\delta$ for all $\delta\in\{0,...,\diam(G)\}$.
The sequence $u_0,...,u_{\diam(G)}$ is called \emph{cosine sequence} of $v$.
\end{observation}

Clearly the cosine sequence and cosine vector of a distance-transitive realization determine each other. We show the following:


\begin{lemma}\label{res:cosine_vectors_are_determined}
\label{res:cosine_vector_of_distance_transitive_realization}
The cosine sequence of a $\theta$-balanced distance-transitive realization (of radius $r(v)=1$) does only depend on $G$ and the eigenvalue $\theta$.
\begin{proof}
Let $N_\delta(i):=\{j\in V\mid\dist(i,j)=\delta\}$ denote the set of all vertices at distance $\delta$ from $i$.
In a distance-transitive graph, the cardinality of the intersection $N_{\delta_1}(i)\cap N_{\delta_2}(j)$ does only depend on $\delta_1,\delta_2$ and $\dist(i,j)$.
The following parameters are therefore well-defined:
%
$$
c_\delta := |\underbrace{N_{\delta-1}(i)\cap N_1(j)}_{=:N_c}|,\quad 
a_\delta := |\underbrace{N_\delta(i)\cap N_1(j)}_{=:N_a}|,\quad 
b_\delta := |\underbrace{N_{\delta+1}(i)\cap N_1(j)}_{=:N_b}|,
$$
%
whenever $\dist(i,j)=\delta$\footnote{The order of the parameter names $a$, $b$ and $c$ might appear counter intuitive, but is standard in the~litera\-ture (these parameters are used to define distance-regular graphs), and so we shall adopt it here. 
}.
The list of the parameters $a_\delta,b_\delta$ and $c_\delta$ is called the~\emph{inter\-section array} of $G$.
Note that $N(j)=N_a\cupdot N_b\cupdot N_c$.

Now suppose that $v$ is a balanced distance-transitive realization with eigenvalue $\theta\in \Spec(G)$ and cosine sequence $u_\delta$. 
Then for all $\delta\in\{0,...,\diam(G)\}$ there is an $i\in N_\delta(1)$, and from that we derive
\vspace{-0.5em}
\begin{align*}
\theta u_\delta 
=  \<v_1,\theta v_i\> 
\overset{\eqref{eq:balanced}} = \Big\< v_1, \sum_{\mathclap{j\in N(i)}} v_j\Big\> 
&= 
\sum_{\mathclap{j\in N_c}} \overbrace{\overset{\phantom.}\<v_1,v_j\>}^{u_{\delta-1}}+
\sum_{\mathclap{j\in N_a}} \overbrace{\overset{\phantom.}\<v_1,v_j\>}^{u_{\delta}}+
\sum_{\mathclap{j\in N_b}} \overbrace{\overset{\phantom.}\<v_1,v_j\>}^{u_{\delta+1}}
\\ &= c_\delta u_{\delta-1} + a_\delta u_\delta + b_\delta u_{\delta+1}.
\end{align*}
Rearranging for $u_{\delta+1}$ yields a three term recurrence for the components of the cosine sequence that only involves $\theta$ and the intersection array:
\begin{equation}
\label{eq:recurrence}
u_{\delta+1} = \frac1{b_\delta}\big((\theta-a_\delta)u_\delta-c_\delta u_{\delta-1}\big).
\end{equation}
We assume that $r(v)=1$, and since $v$ is also arc-transitive we have initial conditions:
$$u_0=[r(v)]^2 = 1,
\qquad 
u_1=\omega(v)\overset{\eqref{eq:relative}}=\frac{\theta}{\deg(G)}.$$
The initial conditions only depends on $\theta$ and $G$, and so the whole cosine sequence does only depend on $\theta$ and $G$ (its degree, and intersection array).
\end{proof}
\end{lemma}

\Cref{res:cosine_vectors_are_determined} together with \cref{res:cosine_vectors} shows that all eigenspaces of a distance-transitive graph are $\Sigma$-irreducible for all distance-transitive $\Sigma\subseteq\Aut(G)$.

We close with an example that demonstrates the potential of the cosine vector approach by applying it outside the realm of distance-transitive realizations.

\begin{example}\label{ex:24_cell_cuboctahedron}
We show that the skeleton of the 24-cell is a spectral realization without computing its spectrum (as we have done in \cref{ex:24_cell}).
Note that the 24-cell is arc-transitive, but not distance-transitive, and that its skeleton is balanced by \cref{res:full_local_dimension_is_good}. 

The coordinates of the vertices of the 24-cell are all coordinate permutations and sign selections of
$$\sqrt2\cdot (\pm 1,0,0,0)\wideand \sqrt2\cdot (\pm1/2,\pm1/2,\pm1/2,\pm1/2).$$%
From this we find that the cosine vector is $u=(2^1,1^8,0^6,(-1)^8,(-2)^1)$ (ignoring the exact ordering of the entries, only caring about the multiplicities).
%

Note that the single entry of value 2 in the cosine vector belongs to the radius $r(v)=\sqrt 2$ of this realization. Also, the eight entries with value $\<v_1,v_i\>=1$ belong to the neighbors $i\in N(1)$, and so this value is determined by \eqref{eq:relative}.
In conclusion, any other balanced arc-transitive realization $\bar v$ to the same eigenvalue (and of the same radius) must have a cosine vector of the form
$$\bar u=(2^1,1^8;x_1,...,x_6;y_1,...,y_8;z),$$
where the $x_i$ match up with the 0-entries in $u$, the $y_i$ match up with the $-1$-entries in $u$, and $z$ matches with the $-2$-entry.

%
As discussed in \cref{res:cosine_vectors}, we can assume $\<u,\bar u\>=0$, which expands to
$$(*)\quad 0=\<u,\bar u\> = 4 + 8  - y_1-\cdots -y_8 - 2z.$$
A (full-dimensional) arc-transitive realization is always centered at the origin, which means $v_1+\cdots +v_n=0$, or in terms of the cosine vector
$$(**)\quad 0=\sum_{i\in V} \bar u_i = 2+8 + x_1+\cdots + x_6 + y_1+\cdots +y_8 + z $$
We can add $(*)$ and $(**)$ to obtain $z=22+x_1+\cdots+ x_6$.

Finally, every component $\bar u_i=\<\bar v_1,\bar v_i\>$ of the cosine vectors must satisfy $-2=\<\bar v_1,-\bar v_1\>\le \< \bar v_1,\bar v_i\> \le \<\bar v_1,\bar v_1\>=2$, thus $x_i,z\in[-2,2]$.
But this is incompatible with $z=22+x_1+\cdots+ x_6$.
%
Thus, no second such realization can exist, and the skeleton of the 24-cell is a spectral realization.

Note that we essentially used the shape of the cosine vector of the 24-cell.
The same argument works essentially unchanged \eg\ for the skeleton of the cuboctahedron (also arc-transitive) whose cosine vector is $(2^1,1^4,0^2,(-1)^4,(-2)^1)$,

\vspace{0.5em}
\begin{center}
\includegraphics[width=0.2\textwidth]{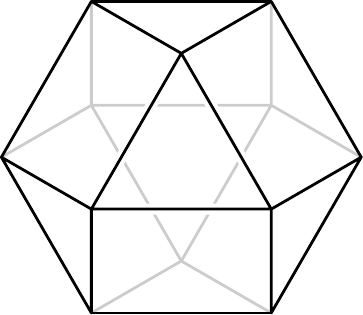}
\end{center}
\end{example}

%% file: sec/future.tex
\tempnewpage
\section{Conclusions and open questions}
\label{sec:future}


In this paper we have taken a look at certain classes of highly symmetric graphs, their symmetric, balanced and spectral realizations.
We were able to show that certain sufficiently symmetric realizations are necessarily spectral, or at least balanced.

We list some open questions.
Most of these questions have a formulation in the language of realizations, and another one in the language of arrangement spaces, and can therefore be attacked from both directions.



\begin{question}\label{q:arc_transitive_balanced}
Is every irreducible arc-transitive realization balanced?

Equivalently, are the irreducible subspaces of an arc-transitive group $\Sigma\subseteq\Aut(G)$ always contained in the eigenspaces of $G$?
\end{question}

We have seen that the answer is \emph{yes} if the realization is rigid (\cref{res:rigid_implies_balanced}), full-dimen\-sional (\cref{res:full_local_dimension_is_good}), or distance-transitive (\cref{res:distance_transitive_consequences}).

\begin{question}\label{q:arc_transitive_rigid}
Are arc-transitive $\Aut(G)$-realizations always rigid?

Equivalently, if $G$ is arc-transitive, are there only finitely many $\Aut(G)$-invariant subspaces of $\RR^n$?
\end{question}

The answer is \emph{yes}, if the realization is distance-transitive (\cref{res:distance_transitive_consequences}).

Another interesting question was asked by Du and Fan in \cite{du2013graph}.
It is known that every group appears as the symmetry group of a graph \cite{frucht1939herstellung}, and with only a few exceptions, almost every group appear as the symmetry group of a \emph{vertex-transitive} graph \cite{babai1995automorphism}.
The eigenspaces of such a graph provide a first clue for the placement of the irreducible subspaces of that initial group, but as we have seen, these connections are not necessarily one-to-one.

\begin{question}
For which groups $\Sigma\subseteq\Sym(n)$ can we find a graph whose eigenspaces are exactly the $\Sigma$-irreducible subspaces of $\RR^n$?
\end{question}

We can certainly do this if $\Sigma=\Aut(G)$ for some distance-transitive graph $G$ as seen in \cref{res:eigenspaces_eq_irreducible}.
If such a graph can be constructed efficiently, this would provide an effective tool for computing invariant subspaces of permutation groups.


\begin{question}\label{q:dim_fix_ge_2}
Can we classify the realizations that are rigid but have $\dim\Fix(T,\Sigma_i)$ $\ge 2$ for all $i\in V$ (\cf\ \cref{res:fix_rigid}).
\end{question}

The classification of these is linked to the exceptional group structures on spheres that only exist on spheres of dimensions $d\in\{0,1,3\}$.
For example,  the realization in \cref{ex:exception} corresponds to the case $d=1$.
The exact connection is not clear to the author.

\subsection{The orbital technique}
\label{sec:orbitals}

In \cref{sec:realizations} we mentioned that symmetric realizations might not be spectral, but only balanced, and that this is a factor preventing us from obtaining all symemtric realizations with spectral methods alone.

We now briefly describe a technique that might fix this problem, though~we~only have empirical evidence for that.
For example, this technique was used to find the balanced, but non-spectral realizations of the graphs in \cref{ex:hexagonal_prism} and \cref{ex:C6_C6}.

The symmetry group $\Aut(G)$ of a graph acts element-wise on the sets $\{i,j\}\subseteq V$ with $i,j\in V$ (note that $\{i,j\}$ can be a singleton if $i=j$).
An \emph{orbital}\footnote{This is non-standard. Usually, an \emph{orbital} is an orbit of $G$ acting on $V\x V$ rather than ${V\choose 2}\cup V$.} is an orbit of this action, that is, it is of the form
$$\mathrm{orb}\{i,j\}:=\{\{\sigma(i),\sigma(j)\}\mid \sigma\in\Aut(G)\}.$$
Let $\mathcal O$ denote the set of all orbitals of $\Aut(G)$.
Consider some map $\mathcal O\ni o\mapsto x_o\in \RR$ that assigns a real number to each  orbital (randomly chosen, or, say, algebraically independent) and define the \emph{orbital matrix} $A^{\mathcal O}\in\Bbb R^{n\x n}$ with entries
$$A^{\mathcal O}_{ij} = x_{\mathrm{orb}\{i,j\}}.$$
Numerical experiments suggest, that the eigenspaces of this matrix are exactly the irreducible invariant subspaces of $\RR^n$ \wrt\ $\Aut(G)$.

While this technique is not perfect (it requires us to obtain the orbitals, which might be a computationally intensive task), it at least provides a finer decomposition of $\RR^n$ than the eigenspaces of the adjacency matrix $A$.

\begin{question}
Are the eigenspaces of the orbital matrix exactly the $\Aut(G)$-irredu\-cible subspaces of $\RR^n$?
\end{question}

%% file: sec/appendix.tex
\appendix

\section{Representation theory}
\label{sec:appendix_rep_theory}

Let $\Sigma\subseteq\Sym(V)$ be a permutation group on $V:=\{1,...,n\}$.

\begin{definition}
A (linear, orthogonal) \emph{$\Sigma$-representation} (or just \emph{representation}) is a group homomorphism $T\:\Sigma\to\Ortho(\RR^d)$, that is
$$T_{\sigma\circ\rho} = T_\sigma T_\rho,\quad\text{for all $\sigma,\rho\in\Sigma$}.$$
\end{definition}

\begin{definition}
Given a representation $T\:\Sigma\to\Ortho(\RR^d)$.
\begin{myenumerate}
	\item A subspace $U\subseteq\RR^d$ is called \emph{$T$-invariant} (or just \emph{invariant}) if $T_\sigma U=U$ for all $\sigma\in\Sigma$. Note that $\{0\}$ and $\RR^d$ are always invariant subspaces.
	\item An invariant subspace $U\subseteq\RR^d$ is called \emph{$T$-irreducible} (or just \emph{irreducible}) if  $U$ and $\{0\}$ are its only invariant subspaces, otherwise it is called \emph{reducible}.
	\item The representation $T$ is called \emph{irreducible} if $\RR^d$ and $\{0\}$ are its only invariant subspaces (that is, $\RR^d$ is irreducible as $T$-invariant subspace), it is called \emph{reducible} otherwise.
\end{myenumerate}
\end{definition}

\begin{remark}\label{rem:intersection}
The intersection of two invariant subspaces $U,U'\subseteq\RR^d$ is again~an invariant~sub\-space, in particular, it is  a subspace of both $U$ and $U'$.
Consequently, if $U$ is irreducible, then either $U\subseteq U'$ or $U\cap U'=\{0\}$.
\end{remark}

\begin{remark}\label{rem:complement}
For a $T$-invariant subspace $U\subseteq\RR^n$, its orthogonal complement $U^\bot$ is again a $T$-invariant subspace.
Applied recursively, we find that $\RR^d$ decomposes~as a direct sum
$$\RR^d=U_1\oplus\cdots \oplus U_m$$
of pairwise orthogonal $T$-irreducible subspaces $U_1,..., U_m\subseteq\RR^d$, though this~decom\-position might not be unique.
\end{remark}

\begin{remark}\label{rem:projection}
If $U\subseteq\RR^d$ is a $T$-invariant subspace, then the orthogonal projection $\pi_U$ onto $U$ commutes with $T$: every vectors $x\in\RR^d$ decomposes like $x=u+u'$ with $u\in U$, and $u'\in U^\bot$. Both $U$ and $U^\bot$ are $T$-invariant (see \cref{rem:complement}), and so $T_\sigma x$ decomposes into $T_\sigma u + T_\sigma u'$ with $T_\sigma u\in U$ and $T_\sigma u'\in U^\bot$. Then
\begin{align*}
T_\sigma\pi_U(x) &= T_\sigma \pi_U(u+u') = T_\sigma u \\&= \pi_U(T_\sigma u+ T_\sigma u') = \pi_U(T_\sigma(u+u'))=\pi_U(T_\sigma x).
\end{align*}
for all $\sigma\in\Sigma$.
Consequently, the projection $\pi_U(U')\subseteq U$ of a $T$-invariant subspace $U'\subseteq\RR^d$ onto $U$ is again $T$-invariant. 
If $U$ is irreducible, then this~pro\-jection must be either $U$ or $\{0\}$.
\end{remark}

For every permutation group $\Sigma\subseteq\Sym(n)$, there is a canonical representation $\sigma\mapsto \Pi_\sigma\in\Perm(n)$ on $\RR^n$ by permutation matrices.
A subspace being invariant or irreducible \wrt\ this representation is called \emph{$\Sigma$-invariant} or \emph{$\Sigma$-irreducible} for~short.


\begin{definition}\label{def:equivariant_isomorphic}
Let $T\:\Sigma\to\Ortho(\RR^d)$ and $T'\!\:\Sigma\to\Ortho(\RR^{d'})$ be two representations.
\begin{myenumerate}
	\item A linear map $M\:\RR^d\to\RR^{d'}$ is called \emph{equivariant} (or \emph{interwining map}) \wrt\ to the pair $(T,T')$, if $MT_\sigma=T'_\sigma M$ for all $\sigma\in \Sigma$.
	\item The representations $T$ and $T'$ are called  \emph{isomorphic} if there exists an invertible equivariant map between them.
\end{myenumerate}
\end{definition}


\begin{theorem}[Schur's lemma; real orthogonal version]\quad
\begin{myenumerate}
\item 
Every equivariant map between two irreducible realizations is either the zero map or invertible. In other words, the only equivariant map between non-isomorphic irreducible $\Sigma$-representations is the zero-map.
\item 
If two (orthogonal) representations are isomorphic, then every equivariant map between them is of the form $\alpha X$, where $\alpha\ge 0$ and $X\in\Ortho(\RR^d)$.
\end{myenumerate}
\end{theorem}

If $U\subseteq\RR^d$ is a $T$-invariant subspace, then we can consider the action of $T_\sigma$~on $U$ as a restricted representation $T_U:\Sigma \to \Ortho(U)$.


\begin{theorem}\label{res:non_orthogonal_subspaces}
Suppose that $U,U'\subseteq\RR^d$ are non-orthogonal irreducible $T$-invariant subspaces. Then
\begin{myenumerate}
	\item $\dim(U)=\dim(U')$, and
	\item the restrictions $T_U$ and $T_{U'}$ are isomorphic representation.
\end{myenumerate}
\begin{proof}
Since $U$ and $U'$ are non-orthogonal subspaces, $\pi_U(U')\not=\{0\}$.
But since~they are~irre\-ducible, we must have $\pi_U(U')=U$ (by \cref{rem:projection}), in particular, $\dim(U)\le \dim(U')$.
We can flip $U$ and $U'$ in this argument to obtain $(i)$.

By the preceding arguments, we can consider $\pi_U\: U'\to U$ as an isomorphism between the subspaces $U$ and $U'$.
Since $\pi_U$ commutes with $T_\sigma$ for all $\sigma\in\Sigma$ (see \cref{rem:projection}), it is a non-zero equivariant map between the restrictions $T_U$ and $T_{U'}$.
The representations are then isomorphic by \cref{def:equivariant_isomorphic} $(ii)$, which gives $(ii)$.
\end{proof}
\end{theorem}

\begin{corollary}\label{res:non_orthogonal_projection}
If $U,U'\subseteq\RR^d$ are non-orthogonal $T$-invariant subspaces, and $U'$ is irreducible, then $\bar U:=\pi_U(U')\subseteq U$ is irreducible.
\begin{proof}
Clearly, $\dim(\bar U)\le \dim(U')$.
But $\bar U$ decomposes into irreducible subspaces, one of which, say $\bar U'$, must be non-orthogonal to $U'$, thus satisfies $\dim \bar U'=\dim U'$ by \cref{res:non_orthogonal_subspaces} $(i)$.
So $\bar U=\bar U'$, and $\bar U$ is irreducible.
\end{proof}
\end{corollary}

\section{Spectral graph theory}
\label{sec:appendix_spec_theory}

Let $G=(V,E)$ be a (simple, undirected) graph with vertex set $V=\{1,...,n\}$,~in particular, on $n$ vertices.

In \emph{spectral graph theory}, when referring to eigenvalues, eigenvectors, eigenspaces, or the spectrum of a graph $G$, one actually refers to the respective quantity for some matrix associated with $G$, mostly its \emph{adjacency matrix} $A\in\{0,1\}^{n\x n}$,
$$A_{ij}:=[ij\in E]=\begin{cases}1&\text{if $ij\in E$}\\0&\text{otherwise}\end{cases},$$
or its \emph{Laplacian} $L\in\ZZ^{n\x n}$,
$$L_{ij} = D - A = \begin{cases} \deg(i) & \text{if $i=j$} \\ -1 & \text{if $ij\in E$} \\ \phantom+ 0&\text{otherwise}\end{cases}.$$
where $D$ is the diagonal matrix with $D_{ii}=\deg(i)$ (the degree of the $i$-th vertex).

Usually, the eigenvalues of $A$ are denoted $\theta_1> \cdots >\theta_m$ (in decreasing order), and the eigenvalues of $L$ (the \emph{Laplacian eigenvalues} of $G$) are denoted $\lambda_1 <\cdots < \lambda_m$ (in increasing oder).
It is well known that $\lambda_1=0$, and therefore $L$ is positive~semi-definite.
Furthermore, the multiplicity of $\lambda_1$ agrees with the number of connected components of $G$.

Both matrices are symmetric, and hence their eigenspaces are pairwise orthogonal. For example, in the case of the adjacency matrix, we obtain a decomposition
$$\RR^n = \Eig_G(\theta_1)\oplus\cdots \oplus\Eig_G(\theta_m)$$
of $\RR^n$ into a direct sum of pairwise orthogonal eigenspaces.

In the case that $G$ is a regular graph of degree $\deg(G)$, the definition of the~Lapla\-cian $L$ simplifies to 
$$L:=\deg(G) \Id- A,$$ 
and the eigenvalues are related via $\lambda_i=\deg(G)-\theta_i$ for all $i\in\{1,...,m\}$.
In~parti\-cular, $\theta_1=\deg(G)$, and the~multi\-plicity of $\theta_1$ indicates the number of connected compo\-nent of $G$.
Furthermore, the eigenspace to $\theta_i$ is exactly the eigenspace to $\lambda_i$.
So in the regular case it suffices to study one set of eigenvalues and eigenvectors, as the results translate directly into the other case.

\section{Implementation in Mathematica}
\label{sec:appendix_mathematica}

The following short Mathematica script takes as input a graph $G$ (in the example below, this is the edge-graph of the dodecahedron), and an index $i$ of an eigenvalue.
It then compute the points $v_i$ (\texttt{vert} in the code), \ie\ the vertex-coordinates of the $\theta_i$-spectal realization.
If the dimension turns out to be appropriate, the spectral realization is plotted.

\vspace{0.5em}
\begin{lstlisting}
(* Input: 
  * the graph G, and
  * the index i of an eigenvalue (i = 1 being the largest eigenvalue).
*)
G = GraphData["DodecahedralGraph"];
i = 2;

(* Computation of vertex coordinates 'vert' *)
n = VertexCount[G];
A = AdjacencyMatrix[G];
eval = Tally[Sort@Eigenvalues[A//N], Round[#1-#2,0.00001]==0 &];
d = eval[[-i,2]]; (* dimension of eigenspace *)
vert = Transpose@Orthogonalize@
  NullSpace[eval[[-i,1]] * IdentityMatrix[n] - A];

(* Output: 
  * the graph G, 
  * its eigenvalues with multiplicities, and
  * the spectral realization.
*)
G
Grid[Join[{{$\theta$,"mult"}}, eval], Frame$\to$All]
Which[
  d<2 , Print["Dimension too low, no plot generated."],
  d==2, GraphPlot[G, VertexCoordinates$\to$vert],
  d==3, GraphPlot3D[G, VertexCoordinates$\to$vert,
  d>3 , Print["Dimension too high, 3-dimensional projection is plotted."];
    GraphPlot3D[G, VertexCoordinates$\to$vert[[;;,1;;3]] ]
]
\end{lstlisting}